\documentclass[paper,twocolumn,twoside]{geophysics}
\usepackage{amssymb}
\usepackage{amsmath}
\usepackage{empheq}
\usepackage{amscd}
\usepackage{graphicx}
\usepackage{graphics}
\usepackage[noadjust]{cite}
\usepackage{amsthm}
\usepackage{caption}
\usepackage{multirow}
\usepackage{array}
\usepackage{rotating}
\usepackage{algorithm}
\usepackage{algorithmic}

\usepackage{indentfirst}
\usepackage{tikz}
\usepackage{calc}
\usepackage{epsfig}
\usepackage{natbib}
\usepackage[makeroom]{cancel}
\usepackage{multicol}
\usepackage{csquotes}

\begin{document}
\title{Attenuation imaging by wavefield reconstruction inversion with bound constraints and total variation regularization}

\address{ \footnotemark[1]University of Tehran, Institute of Geophysics, Tehran, Iran, email: h.aghamiry@ut.ac.ir, agholami@ut.ac.ir \\ 
\footnotemark[2]University Cote d'Azur - CNRS - IRD - OCA, Geoazur, Valbonne, France, email: aghamiry@geoazur.unice.fr, operto@geoazur.unice.fr}
\author{Hossein S. Aghamiry \footnotemark[1]\footnotemark[2], Ali Gholami \footnotemark[1] and St\'ephane Operto \footnotemark[2]}

\lefthead{A PREPRINT~~~~~~~~~~~~~~~~~~~~~~~~~~~~~~~~~~~~~~~~~~~~~~~~~~~~~~~~~~~~~~~~~~~~~~~~~~~~~~~~~~Aghamiry et al.}
\righthead{A PREPRINT~~~~~~~~~~~~~~~~~~~~~~~~~~~~~~~~~~~~~~~~~~~~~~~~~~~~~~~~~~~~~Viscoacoustic Wavefield Inversion}




\begin{abstract}
Wavefield reconstruction inversion (WRI) extends the search space of Full Waveform Inversion (FWI) by allowing for wave equation errors during wavefield reconstruction to match the data from the first iteration. Then, the wavespeeds are updated from the wavefields by minimizing the source residuals. Performing these two tasks in alternating mode breaks down the nonlinear FWI as a sequence of two linear subproblems, relaying on the bilinearity of the wave equation. We solve this biconvex optimization with the alternating-direction method of multipliers (ADMM) to cancel out efficiently the data and source residuals in iterations and stabilize the parameter estimation with appropriate regularizations.
Here, we extend WRI to viscoacoustic media for attenuation imaging. Attenuation reconstruction is challenging because of the small imprint of attenuation in the data and the cross-talks with velocities. To address these issues, we recast the multivariate viscoacoustic WRI as a triconvex optimization and update wavefields, squared slowness, and attenuation factor in alternating mode at each WRI iteration. This requires to linearize the attenuation-estimation subproblem via an approximated trilinear viscoacoustic wave equation. The iterative defect correction embedded in ADMM corrects the errors generated by this linearization, while the operator splitting allows us to tailor $\ell{1}$ regularization to each parameter class. A toy numerical example shows that these strategies mitigate cross-talk artifacts and noise from the attenuation reconstruction. A more realistic synthetic example representative of the North Sea validates the method.
\end{abstract}

\section{Introduction}
%
%
Full waveform inversion (FWI) is a high-resolution nonlinear imaging technology which can provide accurate subsurface model by matching observed and calculated waveforms \citep{Tarantola_1984_ISR,Pratt_1998_GNF,Virieux_2009_OFW}. However, it is well acknowledged that it suffers from two main pathologies. The first one is the nonlinearity associated with cycle skipping: When the distance between the observed and calculated data is the least-squares norm of their differences, FWI remains stuck into spurious local minima when the initial velocity model does not allow to match traveltimes with an error lower than half a period. To mitigate cycle skipping, many variants of FWI have been proposed with more convex distances such as those based on matching filters \citep{Warner_2016_AWI,Guasch_2019_AWI} or optimal transport \citep{Metivier_2018_OTM} among others. The second pathology is ill-posedness resulting from uneven subsurface illumination provided by limited-aperture surface acquisitions \citep[e.g., ][]{Tang_2009_TWL} and parameter cross-talks during multiparameter reconstruction \citep[see ][ for a tutorial]{Operto_2013_TLE}. Mitigating this ill-posedness requires to account for the Hessian in local optimization methods \citep[e.g.][]{Metivier_2017_TRU} and regularize the inversion with prior information such as physical bound constraints \citep[e.g.][]{Asnaashari_2013_RSF,Duan_2016_EWT}. \\
%
%
Among the methods proposed to mitigate cycle skipping, wavefield reconstruction inversion (WRI) \citep{VanLeeuwen_2013_MLM,vanLeeuwen_2016_PMP} extends the parameter search space of frequency-domain FWI by processing the wave equation as a soft constraint with a penalty method. The resulting wave equation relaxation allows for data fitting with inaccurate velocity models through the reconstruction of data-assimilated wavefields, namely wavefields satisfying the observation equation relating the wavefields to the observations \citep{Aghamiry_2019_AEW}.
The algorithm then updates the model parameters by least-squares minimization of the wave equation errors (or source residuals) so that the assimilated wavefields explain both the wave equation and the data as well as possible. 
Performing wavefield reconstruction and parameter estimation in an alternating mode \citep{VanLeeuwen_2013_MLM} rather than by variable projection \citep{vanLeeuwen_2016_PMP} recasts WRI as a sequence of two linear subproblems as a result of the bilinearity of the wave equation in wavefield and squared slowness. The reader is also referred to \citet{Aghamiry_2019_AMW} for a more general discussion on the bilinearity of the elastic anisotropic wave equation. 
\citet{Aghamiry_2019_IWR} solved this biconvex problem with the alternating direction method of multipliers (ADMM) \citep{Boyd_2011_DOS}. ADMM is an augmented Lagrangian method which makes use of operator splitting and alternating directions to solve convex separable multi-variate constrained problems. 
The augmented Lagrangian function combines a penalty function and a Lagrangian function \citep[][ Chapter 17]{Nocedal_2006_NO}. The penalty function relaxes the constraints during early iterations as in WRI, while the Lagrangian function progressively corrects the constraint violations via the action of the Lagrange multipliers.
The leverage provided by the Lagrange multipliers guarantees to satisfy the constraints at the convergence point with constant penalty parameters \citep{Aghamiry_2019_IWR}. Accordingly, \citet{Aghamiry_2019_IWR} called their approach iteratively-refined WRI (IR-WRI). Alternatives to satisfy the constraints at the convergence point with penalty methods rely on multiplicative \citep{daSilva_2017_WRI} or discrepancy-based \citep{Fu_2017_DPM} approaches.
\citet{Aghamiry_2019_IBC} implemented bounding constraints and total variation (TV) regularization \citep{Rudin_1992_NTV} in IR-WRI with the split Bregman method \citep{Goldstein_2009_SBM} to improve the imaging of large-contrast media, with however undesirable staircase imprints in smooth regions. To overcome this issue and capture both the blocky and smooth components of the subsurface, \citet{Aghamiry_2019_CRO} combine in IR-WRI Tikhonov and TV regularizations by infimal convolution.  \\
%
%
The objective of this paper is to extend frequency-domain IR-WRI to viscoacoustic media for attenuation imaging. Attenuation reconstruction by FWI raises two potential issues. The first is related to the cross-talks between wavespeed and attenuation. The ambiguity between velocity and attenuation perturbation in least-squares migration has been emphasized by \citet{Mulder_2009_AAS}. Many combination of velocity and attenuation perturbations can fit equally well reflection amplitudes since they are basically related by an Hilbert transform. This ambiguity can be simply illustrated by the radiation pattern of velocity and attenuation perturbations, which have the same amplitude versus angle behavior and a 90$^\circ$ phase shift \citep{Malinowski_2011_HSA,daSilva_2019_SVF}. The conclusions of \citet{Mulder_2009_AAS} are substantiated by \citet{Ribodetti_2000_AVD} who show that the Hessian of ray+Born least-squares migration of single-offset reflection data is singular if the reflector is not illuminated from above and beneath. On the other hand, \citet{Hak_2011_SAI} show that wavespeed and attenuation can be decoupled during nonlinear waveform inversion of multi-offset/multi-frequency data provided that the causality term is properly implemented in the attenuation model. This conclusion has been further supported by several realistic synthetic experiments and real data case studies in marine and land environments, which manage to reconstruct trustworthy attenuation models \citep{Hicks_2001_RWI,Askan_2007_FWI,Malinowski_2011_HSA,Takougang_2012_SVA,Kamei_2008_WTS,Prieux_2013_MFWa,Stopin_2016_AVF,Operto_2018_MFF,Lacasse_2019_AHA}. This decoupling between velocity and attenuation can be further argued on the basis of physical considerations. In the transmission regime of wave propagation, wavespeeds control the kinematic of wave propagation. This implies that FWI is dominantly driven toward wavespeed updating to match the traveltimes of the wide-aperture data (diving waves, post-critical reflections) and update the long wavelengths of the subsurface accordingly, while attenuation has a secondary role to match amplitude and dispersion effects \citep[e.g., see][ for an illustration]{Operto_2018_MFF}. This weak imprint of the attenuation in the seismic response was illustrated by the sensitivity analysis carried out by \citet{Kurzmann_2013_AFW} who concluded that a crude homogeneous background attenuation model may be enough to perform reliable FWI, while \citet{daSilva_2019_SVF} proposed to reconstruct an under-parametrized attenuation model by semi global FWI. When a high-resolution attenuation model is sought, the ill-posedness of the attenuation reconstruction may be managed with different recipes including data-driven and model-driven inversions (joint versus sequential updates of the velocity and attenuation of selected subdatasets), parameter scaling, bound constraints and regularizations \citep[e.g. ][]{Prieux_2013_MFWa,Operto_2013_TLE}. \\
%
%
In this context, the contribution of this study is two fold: first, we show how to implement velocity and attenuation reconstruction in frequency-domain viscoacoustic IR-WRI when equipped with bound constraints and nonsmooth regularizations. Second, we discuss with numerical examples whether the alternating-direction algorithm driven by the need to expand the search space is suitable to manage ill-conditioned multi-parameter reconstruction.  
It is well acknowledged that viscous effects are easily included in the time-harmonic wave equation with frequency-dependent complex-valued velocities as function of phase velocity and attenuation factor (the inverse of quality factor) which are both real-valued parameters \citep{Toksoz_1981_GRS}. Accordingly, the objective function of viscoacoustic IR-WRI requires to be minimized over a set of three parameter classes (wavefield, squared slowness, attenuation factor). In this study, we consider the Kolsky-Futterman model as attenuation model \citep{Kolsky_1956_PSP,Futterman_1962_DBW}. With this model, the viscoacoustic wave equation is bilinear in wavefield and squared slowness, while it is nonlinear in attenuation factor. This prompts us to introduce a first-order approximation of the viscoacoustic function to form a trilinear viscoacoustic wave equation. This equation allows us to recast the multivariate viscoacoustic IR-WRI as a sequence of three linear subproblems for wavefields, squared slowness and attenuation factor estimation, which are solved in alternating mode following the block relaxation strategy of ADMM. Then, the errors generated by the approximated wave equation during attenuation estimation are corrected by the action of the Lagrange multipliers (dual variables), which are formed by the source residuals computed with the exact wave equation.
Another application of augmented Lagrangian method in AVO inversion is presented in \citet{Gholami_2017_CNA} where the linearized Zoeppritz equations are used to simplify the primal problem, while the dual problem compensates the linearization-related errors by computing the residuals with the exact Zoeppritz equations. Also, the decomposition of the viscoacoustic IR-WRI into three linear subproblems provides the suitable framework to tailor $\ell{1}$ regularizations to each parameter-estimation subproblem \citep{Aghamiry_2019_IBC,Aghamiry_2019_CRO}. \\
%
%
The alternating update of the squared slowness and attenuation factor at each IR-WRI iteration is probably non neutral on how the inversion manages the parameter cross-talks and the contrasted sensitivity of the data to each parameter class, as the multi-parameter inversion is broken down as a sequence of two mono-parameter inversions. This approach differs from those commonly used in multi-parameter inversion. The most brute-force approach consists in the joint updating of the multiple parameter classes, with the issue of managing multi-parameter Hessian with suitable parameter scaling \citep[e.g., ][]{Stopin_2014_MWI,Metivier_2015_AMF,Yang_2016_SVD}. Other approaches rely on ad-hoc hierarchical data-driven and model-driven inversion where the dominant parameter is updated during a first mono-parameter inversion, before involving the secondary parameter in a subsequent multi-parameter inversion \citep[e.g., ][]{Prieux_2013_MFWa,Cheng_2016_MEA}. Another possible model-driven strategy consists of performing the joint updating of the multiple parameter classes during a first inversion, then reset the secondary parameters to their initial values and restart a multi-parameter inversion involving all the parameter classes \citep{Yang_2014_MFW}. \\
In this study, we assess our approach against two synthetic experiments, a toy example and a more realistic well-documented synthetic example representative of the North Sea \citep{Prieux_2013_MFWa}. A comparison of our approach with those reviewed above remains however beyond the scope of this paper. One reason is that all the above approaches have been implemented in conventional FWI, which would remain stuck in a local minimum when starting from the crude initial model used in this study. This is to remind that IR-WRI provides a practical framework to conciliate the search space expansion to mitigate cycle skipping and easy-to-design multi-parameter reconstruction via the alternating update of the multiple parameter classes. \\
%
%
This paper is organized as follow. In the method section, we first review the forward problem equation before going into the details of viscoacoustic IR-WRI: We first formulate the constrained optimization problem to be solved and recast it as a saddle point problem with an augmented Lagrangian function (Appendix A). Then, we review the solution of the three primal subproblems for wavefields, squared slowness and attenuation factor in the framework of ADMM, as well as the expression of the dual variables or Lagrange multipliers that capture the history of the solution refinement in iterations. The appendix B reviews in a general setting the split Bregman method to solve $\ell{1}$-regularized convex problem.  This recipe can be easily applied to the squared slowness and attenuation reconstruction subproblems. The final section presents two synthetic examples, which are performed without and with bound constraints and TV regularization in order to discriminate the role of the augmented-Lagrangian optimization from that of the priors. A toy example allows us to illustrate in a simple setting how well IR-WRI manages the parameter cross-talk and the ill-posedness of the attenuation reconstruction and how bound constraints and TV regularization remove the corresponding artifacts. A second synthetic example representative of the North Sea environment allows one to assess the method in a more realistic setting. 
\section{Theory}
In this section, we first review the viscoacoustic wave equation in the frequency-space domain. Then, we use this wave equation to formulate the iteratively-refined wavefield reconstruction inversion (IR-WRI) for velocity and attenuation.
\subsection{Forward problem}
The viscoacoustic wave equation in the frequency-space domain is given by
\begin{equation} \label{PDEc}
\left(\Delta + \frac{\omega^2}{c(\bold{x})^2} \right)u(\bold{x},\omega) = b(\bold{x},\omega),
\end{equation}
where $\Delta$ is the Laplacian operator, $\omega$ is the angular frequency, $\bold{x}=(x,z)$ denotes the position in the subsurface model and $ b(\bold{x},\omega)$ and $u(\bold{x},\omega)$ are respectively the source term and the wavefield for frequency $\omega$. 
Viscoacoustic (attenuative) media can be described by complex-valued velocity $c(\bold{x})$. 
The velocity associated with the Kolsky-Futterman model is given by \citep{Kolsky_1956_PSP,Futterman_1962_DBW} 
\begin{equation}  \label{KF}
\frac{1}{c(\bold{x})}=\frac{1}{v(\bold{x})}\left[ 1- \frac{1}{\pi Q(\bold{x})} \log |\frac{\omega}{\omega_r}| + i \frac{\text{sign}(\omega)}{2Q(\bold{x})}\right],
\end{equation} 
%
%
where $v(\bold{x})$ denotes the phase velocity, $Q(\bold{x})$ the frequency-independent quality factor, both real-valued, and $i=\sqrt{-1}$. Also, $\text{sign}(\bullet)$ is the sign function that extracts the sign of a real number $\bullet$. 
The logarithmic term with reference frequency $\omega_r=2\pi f_r$ implies causality \citep{Aki_2002_QST,Hak_2011_SAI}. In this study, $f_r$ is chosen to be 50 Hz \citep{Toverud_2005_CSA}.

\subsection{Inverse problem}
We discretize the 2D partial-differential equation (PDE), equation \ref{PDEc}, with a $N = N_x \times N_z$ grid points, where $N_x$ and $N_z$ are the number of points in the horizontal and vertical directions, respectively.
We parametrize the inversion by squared slowness $m=1/v^2$ and attenuation factor
$\alpha=1/Q$. Accordingly, equation \ref{KF} in discrete form reads as
\begin{equation} \label{KFmodel}
\frac{1}{\bold{c}^2}=\bold{m} \circ {\rho}(\boldsymbol{\alpha}),
\end{equation} 
where
\begin{equation} \label{rho}
\rho(\alpha) = \left( 1 + \beta(\omega) \alpha\right)^2, ~~ 
\beta(\omega)=i\frac{ \text{sign}(\omega)}{2} - \frac{1}{\pi} \log |\frac{\omega}{\omega_r}|
\end{equation}
and  $\circ$ denotes Hadamard (element wise) product operator.

The model parameters $\bold{m}\in \mathbb{R}^N$ and $\boldsymbol{\alpha}\in \mathbb{R}^N$ are defined as solution of the following nonlinear PDE-constrained optimization problem \citep{Aghamiry_2019_IWR}
\begin{align} 
&\min_{\bold{u},\bold{m}\in \mathcal{M},\boldsymbol{\alpha}\in \mathcal{A}}~~~~\mu E(\bold{m})+\nu F(\boldsymbol{\alpha}) \nonumber \\
&\text{subject to}
~~~~\begin{cases} \bold{A}(\bold{m},\boldsymbol{\alpha})\bold{u}=\bold{b} \\
\bold{Pu}=\bold{d},
\end{cases}
\label{main}
\end{align} 
where $\bold{u} \in \mathbb{C}^{N\times 1}$ is the wavefield, 
$\bold{b} \in \mathbb{C}^{N\times 1}$ is the source term, and $\bold{A}(\bold{m},\boldsymbol{\alpha}) \in \mathbb{C}^{N\times N}$ is the matrix representation of the discretized Helmholtz PDE, equation \ref{PDEc}.
The observation operator $\bold{P} \in \mathbb{R}^{M\times N}$ samples the reconstructed wavefields at the $M$ receiver positions for comparison with the recorded data $\bold{d} \in \mathbb{C}^{M\times 1}$ (we assume a single source experiment for sake of compact notation; However, the extension to multiple sources is straightforward).
The functions $E$ and $F$ are appropriate regularization functions for $\bold{m}$ and $\boldsymbol{\alpha}$, respectively, which are weighted by the penalty parameters $\mu$ and $\nu>0$, respectively. 
$\mathcal{M}$ and $\mathcal{A}$ are convex sets defined according to our prior knowledge of $\bold{m}$ and $\boldsymbol{\alpha}$. For example, if we know the lower and upper bounds on $\bold{m}$ and $\boldsymbol{\alpha}$ then
\begin{equation}
\mathcal{M} = \{\bold{m} \vert \bold{m}_{min} \leq \bold{m} \leq \bold{m}_{max}\},
\end{equation}
and
\begin{equation}
\mathcal{A} = \{\boldsymbol{\alpha} \vert \boldsymbol{\alpha}_{min} \leq \boldsymbol{\alpha} \leq \boldsymbol{\alpha}_{max}\}.
\end{equation}

The PDE constraint $\bold{A}(\bold{m},\boldsymbol{\alpha})\bold{u}=\bold{b}$ in equation \ref{main} is nonlinear in $\bold{m}$ and $\boldsymbol{\alpha}$ and very ill-conditioned \citep{Dolean_2015_IDD}, while the data constraint $\bold{Pu}=\bold{d}$ is linear but the operator $\bold{P}$ is rank-deficient with a huge null space because $M\ll N$. Therefore, determination of the optimum multivariate solution ($\bold{u}^*,\bold{m}^*,\boldsymbol{\alpha}^*$) satisfying both constraints (the wave equation and the observation equation) simultaneously is extremely difficult, and requires sophisticated regularizations. In this paper, we use the first-order isotropic TV regularization \citep{Rudin_1992_NTV} for both $\bold{m}$ and $\boldsymbol{\alpha}$, i.e. $E(\bold{m})=\|\bold{m}\|_{\text{TV}}$ and $F(\boldsymbol{\alpha})=\|\boldsymbol{\alpha}\|_{\text{TV}}$. However, other regularizations such as compound regularizations can be used in a similar way \citep[see][]{Aghamiry_2019_CRO}.
The isotropic TV norm of a 2D image $\bold{w}\in \mathbb{R}^N$ is defined as \citep{Rudin_1992_NTV}
\begin{equation}
\| \bold{w} \|_\text{TV}=\sum  \sqrt{(\nabla_{\!x} \bold{w})^2 + (\nabla_{\!z}\bold{w})^2},
\end{equation}
where $\nabla_{\!x}$ and $\nabla_{\!z}$ are respectively first-order difference operators in the horizontal and vertical directions with appropriate boundary conditions \citep{Gholami_2019_3DD}. 

Beginning with an initial velocity model $\bold{v}^0=1/\sqrt{\bold{m}^{0}}$ and $\boldsymbol{\alpha}^{0}=\bold{0}$, $\bold{b}^{0}=\bold{0}$, $\bold{d}^{0}=\bold{0}$, ADMM solves iteratively the multivariate optimization problem, equation \ref{main}, with alternating directions as \citep[see][ for more details]{Boyd_2011_DOS,Benning_2015_ADM,Aghamiry_2019_IWR,Aghamiry_2019_IBC}
\begin{subequations}
\label{ADMM}
 \begin{empheq}[left={\empheqlbrace\,}]{align}
\bold{u}^{k+1}&= \underset{\bold{u}}{\arg\min} ~ \Psi(\bold{u},\bold{m}^k,\boldsymbol{\alpha}^k,\bold{b}^k,\bold{d}^k) \label{primal_u}\\
\bold{m}^{k+1}&= \underset{\bold{m}\in \mathcal{M}}{\arg\min} ~ \Psi(\bold{u}^{k+1},\bold{m},\boldsymbol{\alpha}^k,\bold{b}^k,\bold{d}^k) \label{primal_sigma}\\
\boldsymbol{\alpha}^{k+1}&= \underset{\boldsymbol{\alpha}\in \mathcal{A}}{\arg\min} ~ \Psi(\bold{u}^{k+1},\bold{m}^{k+1},\boldsymbol{\alpha},\bold{b}^k,\bold{d}^k) \label{primal_eta}\\
\bold{b}^{k+1} &= \bold{b}^k  +\bold{b}- \bold{A}(\bold{m}^{k+1},\boldsymbol{\alpha}^{k+1})\bold{u}^{k+1} \label{dual_b}\\ 
\bold{d}^{k+1} &= \bold{d}^k  +\bold{d}- \bold{P}\bold{u}^{k+1},  \label{dual_d}
\end{empheq}
\end{subequations}
where
\begin{align}
& ~~~~\Psi(\bold{u},\bold{m},\boldsymbol{\alpha},\bold{b}^k,\bold{d}^k) =
\mu \|\bold{m}\|_{\text{TV}}+\nu \|\boldsymbol{\alpha}\|_{\text{TV}} \label{eqpsi}\\
&+ \lambda \|\bold{b}^k+\bold{b}-\bold{A}(\bold{m},\boldsymbol{\alpha})\bold{u}\|_2^2 + 
\gamma\|\bold{d}^k+\bold{d}-\bold{Pu}\|_2^2, \nonumber
\end{align}
is the augmented Lagrangian function written in scaled form (Appendix A), $\bullet^k$ is the value of $\bullet$ at iteration $k$, the scalars $\lambda,\gamma>0$ are the penalty parameters assigned to the wave-equation and observation-equation constraints, respectively, and $\bold{b}^{k}$, $\bold{d}^{k}$ are the scaled Lagrange multipliers, which are updated through a dual ascent scheme by the running sum of the constraint violations (source and data residuals) as shown by  equations~\ref{dual_b}-\ref{dual_d}. The penalty parameters $\lambda,\gamma>0$ can be tuned in equation \ref{eqpsi} such that a dominant weight $\gamma$ is given to the observation equation at the expense of the wave equation during the early iterations to guarantee the data fit, while the iterative update of the Lagrange multipliers progressively correct the errors introduced by these penalizations such that both of the observation equation and the wave equation are satisfied at the convergence point with acceptable accuracies.
In the next three subsections, we show how to solve each optimization subproblem \ref{primal_u}-\ref{primal_eta}.
\subsubsection{Update wavefield (subproblem \ref{primal_u})}
The objective function $\Psi$ is quadratic in $\bold{u}$ and its minimization gives the following closed-form expression of  $\bold{u}$
\begin{align}
\left(\lambda \bold{A}^T \bold{A} +\gamma \bold{P}^T \bold{P} \right) \bold{u}^{k+1} = 
 \lambda \bold{A}^T (\bold{b}^k+\bold{b}) + \gamma \bold{P}^T (\bold{d}^k+\bold{d}),
\label{eqnormal}
\end{align}
where $\bold{A}\equiv \bold{A}(\bold{m}^k,\boldsymbol{\alpha}^k)$ and $\bold{A}^T$ denotes the Hermitian transpose of $\bold{A}$. 
%

\subsubsection{Update squared slowness (subproblem \ref{primal_sigma})}
The PDE operator
\begin{equation}   \label{Adef} 
\bold{A}(\bold{m},\boldsymbol{\alpha}) =  \bold{\Delta} + \omega^2 \bold{C} \text{diag}(\bold{m}\circ {\rho}(\boldsymbol{\alpha}))\bold{B}
\end{equation}
is discretized with the finite-difference method of \citet{Chen_2013_OFD} where  $\bold{\Delta}$ is the discretized Laplace operator, $\bold{C}$ introduces boundary conditions such as perfectly-matched layers \citep{Berenger_1994_PML}, $\bold{B}$ is the mass matrix \citep{Marfurt_1984_AFF} which spreads the mass term $\omega^2\bold{C} \text{diag}(\bold{m}\circ {\rho}(\boldsymbol{\alpha}))$ 
over all the  coefficients of the stencil to improve its accuracy following an anti-lumped mass strategy,  and
diag($\bullet$) denotes a diagonal matrix.
From equation \ref{Adef} we get that
\begin{equation}   \label{Au}
\bold{A}(\bold{m},\boldsymbol{\alpha})\bold{u} =  \bold{A}(\boldsymbol{0},\boldsymbol{\alpha})\bold{u} + \omega^2 \bold{C} \text{diag}(\bold{Bu}\circ {\rho}(\boldsymbol{\alpha}))\bold{m},
\end{equation}
where $\bold{A}(\boldsymbol{0},\boldsymbol{\alpha})\equiv \bold{\Delta}$. 
Therefore, subproblem \ref{primal_sigma} can be written as
\begin{equation} \label{TV_sig}
\bold{m}^{k+1}= \underset{\bold{m}\in \mathcal{M}}{\arg\min}~
\mu \|\bold{m}\|_{\text{TV}} + \lambda \|\bold{L}\bold{m}-\bold{y}^k\|_2^2,
\end{equation}
where
\begin{equation}
\bold{L} =  \omega^2 \bold{C} \text{diag}(\bold{Bu}^{k+1}\circ {\rho}(\boldsymbol{\alpha}^k)),
\end{equation}
and
\begin{equation}
\bold{y}^k=\bold{b}^k+\bold{b}- \bold{\Delta} \bold{u}^{k+1}.
\end{equation}
Equation \ref{TV_sig} describes the box-constrained TV-regularized subproblem for $\bold{m}$ which is convex but non-smooth. This box-constrained TV-regularized problem can be solved efficiently with ADMM and splitting methods, also referred to as the split Bregman method \citep{Goldstein_2009_SBM}. 
Using splitting methods, the unconstrained subproblem \ref{TV_sig} is recast as a multivariate constrained problem, through the introduction of auxiliary variables. These auxiliary variables are introduced to decouple the $\ell{2}$ subproblem from the $\ell{1}$ subproblem such that they can be solved in  alternating mode with ADMM. Moreover, a closed form expression of the auxiliary variables is easily obtained by solving the $\ell{1}$ subproblem with proximity operators \citep{Combettes_2011_PRO,Parikh_2013_PA}. We refer the reader to Appendix B for a more detailed review of this method.
\subsubsection{Update attenuation factor (subproblem \ref{primal_eta})}
Subproblem \ref{primal_eta} is nonlinear due to the nonlinearity of the PDE with respect to $\boldsymbol{\alpha}$.
We linearize this subproblem by using a first-order approximation of $\rho(\alpha)$, equation \ref{rho}, as
\begin{equation}
\rho(\alpha) \approx 1 + 2\beta(\omega)\alpha,
\end{equation}
which is accurate for $\alpha \ll 1$ \citep{Hak_2011_SAI}, and gives
\begin{equation}   \label{Aua}
\bold{A}(\bold{m},\boldsymbol{\alpha})\bold{u} \approx  \bold{A}(\bold{m},\boldsymbol{0})\bold{u}
 + 2\omega^2\beta(\omega) \bold{C} \text{diag}(\bold{Bu} \circ \bold{m})\boldsymbol{\alpha}.
\end{equation}
Accordingly, subproblem \ref{primal_eta} can be written as the following linear problem
\begin{equation} \label{TV_eta}
\boldsymbol{\alpha}^{k+1}\approx \underset{\boldsymbol{\alpha}\in \mathcal{A}}{\arg\min}~
\nu \|\boldsymbol{\alpha}\|_\text{TV} + \lambda \|\bold{H}\boldsymbol{\alpha}-\bold{h}^k\|_2^2,
\end{equation}
where
\begin{equation}
\bold{H} =  2\omega^2\beta(\omega) \bold{C} \text{diag}( \bold{Bu}^{k+1}\circ\bold{m}^{k+1}),
\end{equation}
and
\begin{equation}
\bold{h}^k=\bold{b}^k+\bold{b}- \bold{A}(\bold{m}^{k+1},\boldsymbol{0})\bold{u}^{k+1}.
\end{equation}

Equation \ref{TV_eta} is also a box-constrained TV-regularized convex problem, which can be solved with ADMM (Appendix A) in a manner similar to the previous subproblem for squared slowness.
It is important to stress that the errors generated by the first-order approximation of $\rho(\boldsymbol{\alpha})$ during the update of $\boldsymbol{\alpha}$ are iteratively compensated by the action of the scaled Lagrange multiplier $\bold{b}^k$.
These Lagrange multipliers are formed by the running sum of the wave equation errors, which are computed with the exact wave equation operator (namely, without linearization of $\rho(\boldsymbol{\alpha})$). \\
The overall workflow described above is summarized in Algorithm 1. 
%

\begin{algorithm*}[htb]
\caption{
Viscoacoustic Wavefield Inversion with Bound Constrained TV Regularization. Lines 4 to 6 are the primal subproblems for wavefield reconstruction and parameter estimation. Lines 7 to 12 are primal subproblems for auxiliary variables introduced to implement nonsmooth regularizations and bound constraints (Appendix B).  Lines 13 to 20 are the dual subproblems solved with gradient ascent steps.}
\label{Alg2cont0}
\footnotesize
{\fontsize{10}{18}\selectfont
\begin{algorithmic}[1]
\STATE Begin with $k=0$, an initial squared slowness $\bold{m}^0$, and attenuation $\boldsymbol{\alpha}^{0}$,
\STATE Set to zero the values of $\bold{d}^0$, $\bold{b}^0,\bold{p}^0_{x,m},\bold{p}^0_{y,m},\bold{p}^0_{z,m},\bold{p}^0_{x,\alpha}, \bold{p}^0_{y,\alpha},\bold{p}^0_{z,\alpha},\bold{q}^0_{x,m},\bold{q}^0_{y,m},\bold{q}^0_{z,m},\bold{q}^0_{x,\alpha}, \bold{q}^0_{y,\alpha},\bold{q}^0_{z,\alpha},$
\WHILE {convergence criteria  not satisfied}
\STATE $\bold{u}^{k+1} ~=  \Big[\lambda \bold{A}^T \bold{A} +\gamma \bold{P}^T \bold{P} \Big]^{-1} 
\Big[ \lambda \bold{A}^T [\bold{b}^k+\bold{b}] + \gamma \bold{P}^T [\bold{d}^k+\bold{d}]\Big]$
\STATE $\bold{m}^{k+1} =\Big[\lambda\bold{L}^T \bold{L} + \xi_m \nabla_{\!x}^T\nabla_{\!x}+\xi_m \bold{I} +\xi_m \nabla_{\!z}^T\nabla_{\!z} \Big]^{-1}\Big[\lambda\bold{L}^T \bold{y}^k + \xi_m \nabla_{\!x}^T [\bold{p}_{x,m}^k+\bold{q}_{x,m}^k]+\xi_m [\bold{p}_{y,m}^k+\bold{q}_{y,m}^k]+\xi_m \nabla_{\!z}^T [\bold{p}_{z,m}^k+\bold{q}_{z,m}^k]\Big]$
\STATE $\boldsymbol{\alpha}^{k+1} = \Big[\lambda\bold{H}^T \bold{H} + \xi_\alpha \nabla_{\!x}^T\nabla_{\!x}+\xi_\alpha \bold{I} +\xi_\alpha \nabla_{\!z}^T\nabla_{\!z} \Big]^{-1}\Big[\lambda\bold{H}^T \bold{h}^k + \xi_\alpha \nabla_{\!x}^T [\bold{p}_{x,\alpha}^k+\bold{q}_{x,\alpha}^k]+\xi_\alpha [\bold{p}_{y,\alpha}^k+\bold{q}_{y,\alpha}^k]+\xi_\alpha \nabla_{\!z}^T [\bold{p}_{z,\alpha}^k+\bold{q}_{z,\alpha}^k]\Big]$
\STATE $\bold{p}_{x,m}^{k+1} = \max(1 - \frac{\mu/\xi_m}{\sqrt{|\nabla_{\!x} \bold{m}^{k+1}-\bold{q}_{x,m}^k|^2+|\nabla_{\!z} \bold{m}^{k+1}-\bold{q}_{z,m}^k|^2}},0)\circ (\nabla_{\!x} \bold{m}^{k+1}-\bold{q}_{x,m}^k)$
\STATE $\bold{p}_{y,m}^{k+1} = \text{proj}_{\mathcal{M}} (\bold{m}^{k+1} - \bold{q}_{y,m}^k)$
\STATE $\bold{p}_{z,m}^{k+1} = \max(1 - \frac{\mu/\xi_m}{\sqrt{|\nabla_{\!x} \bold{m}^{k+1}-\bold{q}_{x,m}^k|^2+|\nabla_{\!z} \bold{m}^{k+1}-\bold{q}_{z,m}^k|^2}},0)\circ (\nabla_{\!z} \bold{m}^{k+1}-\bold{q}_{z,m}^k)$
\STATE $\bold{p}_{x,\alpha}^{k+1} = \max(1 - \frac{\nu/\xi_\alpha}{\sqrt{|\nabla_{\!x} \boldsymbol{\alpha}^{k+1}-\bold{q}_{x,\alpha}^k|^2+|\nabla_{\!z} \boldsymbol{\alpha}^{k+1}-\bold{q}_{z,\alpha}^k|^2}},0)\circ (\nabla_{\!x} \boldsymbol{\alpha}^{k+1}-\bold{q}_{x,\alpha}^k)$
\STATE $\bold{p}_{y,\alpha}^{k+1} = \text{proj}_{\mathcal{A}} (\boldsymbol{\alpha}^{k+1} -\bold{q}_{y,\alpha}^k)$
\STATE $\bold{p}_{z,\alpha}^{k+1} = \max(1 - \frac{\nu/\xi_\alpha}{\sqrt{|\nabla_{\!x} \boldsymbol{\alpha}^{k+1}-\bold{q}_{x,\alpha}^k|^2+|\nabla_{\!z} \boldsymbol{\alpha}^{k+1}-\bold{q}_{z,\alpha}^k|^2}},0)\circ (\nabla_{\!z} \boldsymbol{\alpha}^{k+1}-\bold{q}_{z,\alpha}^k)$
\STATE $\bold{q}_{x,m}^{k+1} = \bold{q}^k_{x,m} + \bold{p}_{x,m}^{k+1}-\nabla_{\!x} \bold{m}^{k+1}$
\STATE $\bold{q}_{y,m}^{k+1} = \bold{q}_{y,m}^k + \bold{p}_{y,m}^{k+1}-\bold{m}^{k+1}$
\STATE $\bold{q}_{z,m}^{k+1} = \bold{q}^k_{z,m} + \bold{p}_{z,m}^{k+1}-\nabla_{\!z} \bold{m}^{k+1}$
\STATE $\bold{q}_{x,\alpha}^{k+1} = \bold{q}^k_{x,\alpha} + \bold{p}_{x,\alpha}^{k+1}-\nabla_{\!x} \boldsymbol{\alpha}^{k+1}$
\STATE $\bold{q}_{y,\alpha}^{k+1} = \bold{q}_{y,\alpha}^k + \bold{p}_{y,\alpha}^{k+1}-\boldsymbol{\alpha}^{k+1}$
\STATE $\bold{q}_{z,\alpha}^{k+1} = \bold{q}^k_{z,\alpha} + \bold{p}_{z,\alpha}^{k+1}-\nabla_{\!z} \boldsymbol{\alpha}^{k+1}$
\STATE $\bold{b}^{k+1} = \bold{b}^k  +\bold{b}- \bold{A}(\bold{m}^{k+1},\boldsymbol{\alpha}^{k+1})\bold{u}^{k+1}$
\STATE $\bold{d}^{k+1} = \bold{d}^k  +\bold{d}- \bold{P}\bold{u}^{k+1}$
\STATE $k = k+1$ 
\ENDWHILE
 \end{algorithmic}
}
\end{algorithm*}
\subsection{Practical implementation}
The ADMM optimization that is used to solve equations \ref{TV_sig} and \ref{TV_eta} is reviewed in Appendix B. We also refer the reader to \citet{Goldstein_2009_SBM}, \citet{Boyd_2011_DOS} and \citet{Aghamiry_2019_IBC} as a complement. A key property of the ADMM algorithm, equation \ref{ADMM}, is that, at iteration $k$, we don't need to solve each optimization subproblems \ref{primal_u}-\ref{primal_eta} exactly via inner iterations. The intuitive reason is that the updating of the primal variable performed by one subproblem is hampered by the errors of the other primal variables that are kept fixed. In this framework, the errors at each iteration $k$ are more efficiently compensated by the gradient-ascent update of the Lagrange multipliers (dual variable). This statement was corroborated by numerical experiments which showed that one (inner) iteration of each subproblem per ADMM cycle $k$ generates solutions which are accurate enough to guarantee the fastest convergence of the ADMM algorithm \citep{Goldstein_2009_SBM,Boyd_2011_DOS,Aghamiry_2019_IBC,Gholami_2017_CNA,Aghamiry_2019_AMW}. Moreover, this error compensation is more efficient when the dual variables are updated after each primal subproblem \ref{primal_u}-\ref{primal_eta} rather than at the end of an iteration $k$ as indicated in Algorithm 1 for sake of compactness. This variant of ADMM is referred to as the Peaceman-Rachford splitting method \citep{peaceman_1955_PRA,He_2014_SPR} and will be used in the following numerical experiments. The reader is referred to \citet{Aghamiry_2019_IWR} for more details about the improved convergence of the Peaceman-Rachford splitting method compared to ADMM in the framework of IR-WRI. 

We follow the guideline presented in \citet[][ section 3.1]{Aghamiry_2019_IBC} for tuning the penalty parameters. The overall procedure is as follow:  we first set $\mu$=0.6 and $\nu$=0.4 to tune the relative weight of the regularizations of the squared slowness ($\bold{m}$) and attenuation factor ($\boldsymbol{\alpha}$), equation~\ref{eqpsi}. Then, we set the ratios $\mu/\xi_m$ and $\nu/\xi_\alpha$, lines 7, 9, 10, 12 in Algorithm 1, to tune the soft thresholding performed by the TV regularization of $\bold{m}$ and $\boldsymbol{\alpha}$ (subproblems \ref{primal_sigma} and \ref{primal_eta}). We refer the reader to Appendix B for the role of the penalty parameters $\xi_m$ and $\xi_\alpha$ in TV regularization (denoted generically by $\xi$ in equation~\ref{CTV2}). We set $\mu/\xi_m$ and $\nu/\xi_\alpha$ equal to 0.02$\times$max($\bold{r})$, equations~\ref{r}. These values can be refined according to our prior knowledge of the subsurface medium. Then,  we set a constant $\lambda$ to balance the relative weight of the regularization and the wave-equation misfit function during the parameter-estimation subproblems, equations~\ref{TV_sig} and \ref{TV_eta}, and lines 5 and 6 in Algorithm 1. If necessary, $\lambda$ can be increased during iterations to mitigate the imprint of the regularization near the convergence point. Finally, we set $\gamma$ for wavefield reconstruction such that $\lambda / \gamma$ is a small fraction of the highest eigenvalue of the regularized normal operator, equation~\ref{eqnormal} and line 1 in Algorithm 1 \citep{vanLeeuwen_2016_PMP}. This parameter $\gamma$ can be kept constant during iterations. The reader is referred to \citet{Aghamiry_2019_IWR} for an analysis of the sensitivity of IR-WRI to the choice of the weight balancing the role of the observation equation and the wave equation during wavefield reconstruction.

\section{NUMERICAL EXAMPLES}
%
%
\subsection{Simple inclusions test}
\begin{figure}[ht!]
\begin{center}
\includegraphics[width=0.49\textwidth]{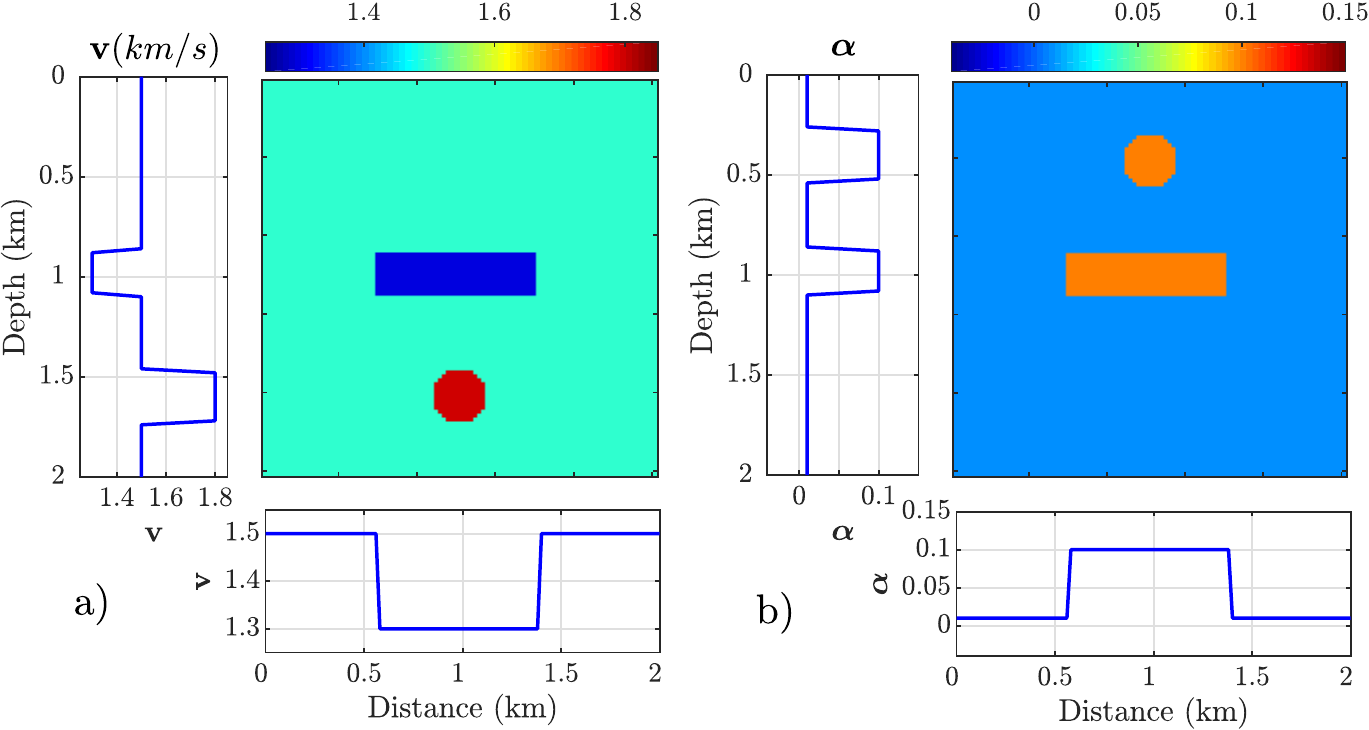}\\
\includegraphics[width=0.49\textwidth]{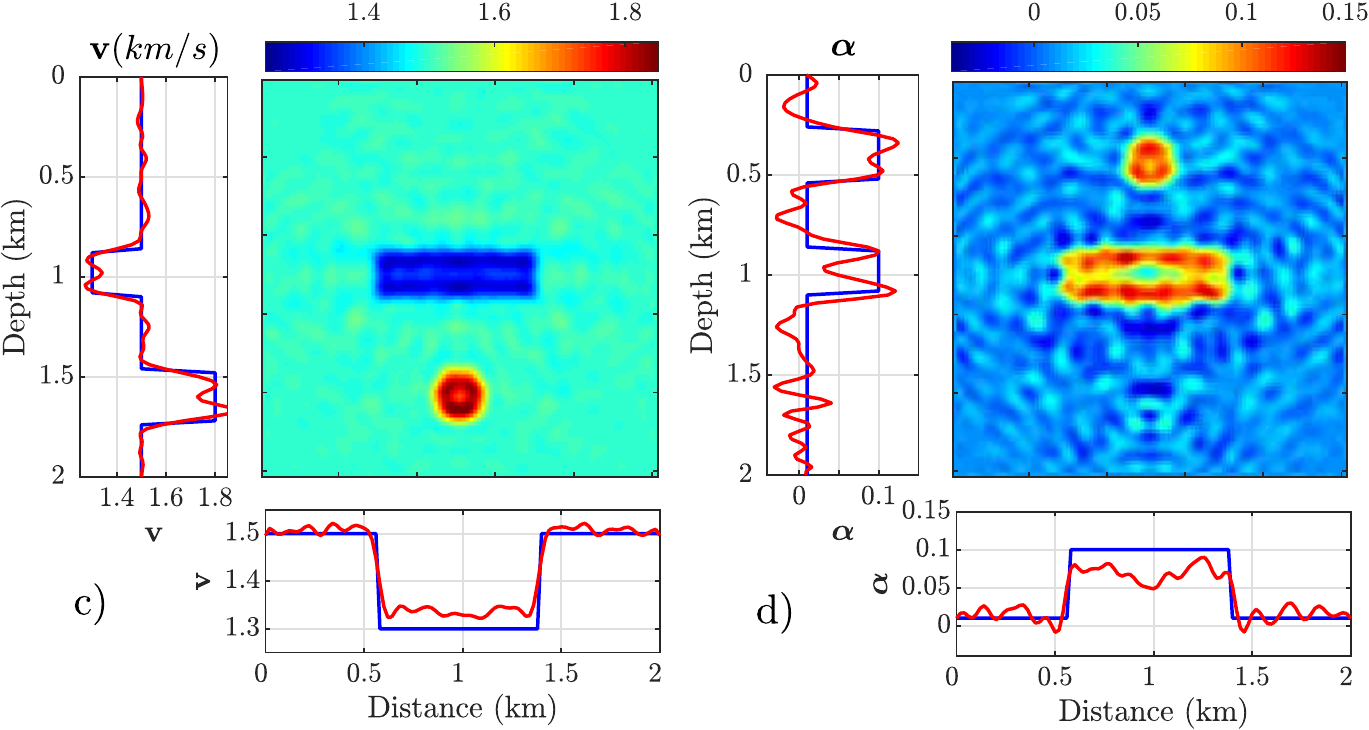}\\
\includegraphics[width=0.49\textwidth]{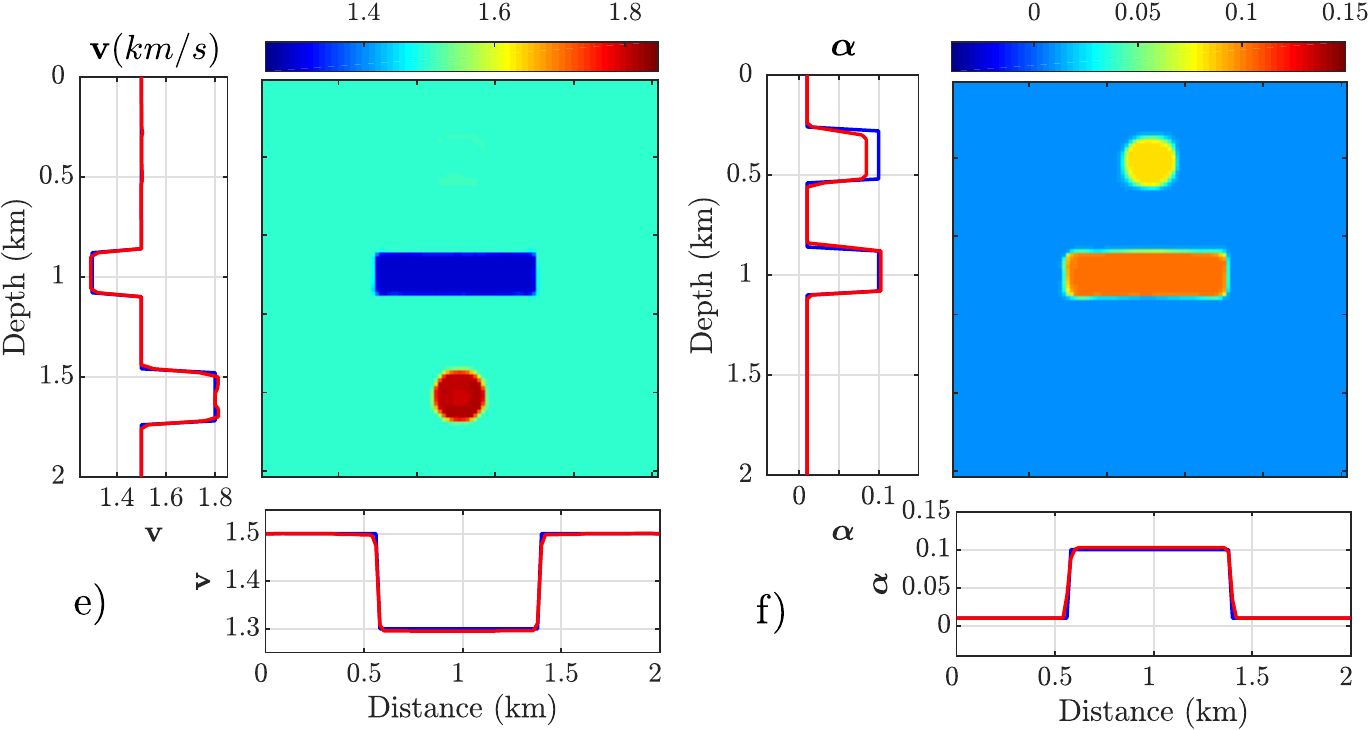}
\caption{(a) True velocity model. (b) True attenuation model. (c-d)  Reconstructed velocity (c) and attenuation (d) from viscoacoustic IR-WRI without TV regularization. (e-f) Same as (a-b) when TV regularization is applied. Profiles of the true (blue) and reconstructed (red) models running across the center of the models are shown on the left and bottom of the reconstructed models.}
\label{fig:Box_test}
\end{center}
\end{figure}

We first consider a simple 2D example to validate viscoacoustic IR-WRI without and with TV regularization. The true velocity model is a homogeneous background model with a wavespeed of 1.5 km/s, which contains two inclusions: a 250-m-diameter circular inclusion at position (1~km,1.6~km) with a wavespeed of 1.8 km/s and a 0.2 $\times$ 0.8 km rectangular inclusion at the center of the model with a wavespeed of 1.3 km/s (Figure \ref{fig:Box_test}a). Also, the true $\boldsymbol{\alpha}$ model is a homogeneous background model with $\boldsymbol{\alpha}$=0.01 ($\bold{Q}=100$), which contains two inclusions, both of them with $\boldsymbol{\alpha}=0.1$ ($\bold{Q}=10$). The first is a 250-m-diameter circular inclusion at position (1~km,0.4~km) and the second is a 0.2 $\times$ 0.8 km rectangle inclusion at the center of the model (Figure \ref{fig:Box_test}b). The wavespeed and $\boldsymbol{\alpha}$ rectangular inclusions share the same size and position, while the position of the wavespeed and $\boldsymbol{\alpha}$ circular inclusions are different in order to test different parameter trade-off scenarios.
A vertical and horizontal logs which cross the center of $\bold{v}$ and $\boldsymbol{\alpha}$ models are plotted in the left and bottom side of the models, respectively, for all the figures of this test. 
An ideal acquisition is used with 8 sources and 200 receivers along the four edges of the model and three frequency components (2.5, 5, 7~Hz) are jointly inverted using a maximum of 30 iterations as stopping criterion of iteration. 

We performed viscoacoustic IR-WRI without and with TV regularization starting from homogeneous $\bold{v}$ and $\boldsymbol{\alpha}$ models with  $\bold{v}$ = 1.5 km/s and $\boldsymbol{\alpha} = {0}$, respectively. The final $(\bold{v},\boldsymbol{\alpha})$ models estimated by IR-WRI without and with TV regularization are shown in Figures~\ref{fig:Box_test}c-d and \ref{fig:Box_test}e-f, respectively. 
%
%
The IR-WRI results without TV regularization show acceptable velocity model and quite noisy attenuation reconstruction (Figure~\ref{fig:Box_test}c-d). The velocity reconstruction is however hampered by significant limited bandwidth effects. We also show underestimated wavespeeds in the rectangular inclusion along the horizontal log. These underestimated wavespeeds are clearly correlated with underestimated $\boldsymbol{\alpha}$ values (overestimated $\bold{Q}$), hence highlighting some trade-off between $\bold{v}$ and $\boldsymbol{\alpha}$. Other moderate trade-off artifacts are shown in the vertical log of the $\boldsymbol{\alpha}$ model, which shows undesired high-frequency perturbations at the position of the circular velocity inclusion.
The TV regularization removes very efficiently all of these pathologies: it extends the wavenumber bandwidth of the models and removes to a large extent the parameter cross-talk as the subsurface medium perfectly matches the  piecewise-constant prior associated with TV regularization (Figure \ref{fig:Box_test}e-f). 
Note however that $\boldsymbol{\alpha}$ remains slightly underestimated in the circular inclusion (Figure~\ref{fig:Box_test}f). This correlates with a barely-visible velocity underestimation at this location in Figure~\ref{fig:Box_test}e. The relative magnitude of these errors gives some insight on the relative sensitivity of the data to $\bold{m}$ and $\boldsymbol{\alpha}$.

%
%
%
%
%
%
%
%
\subsection{Synthetic North Sea case study}
\subsubsection{Experimental setup}
\begin{figure}[ht!]
\centering
\includegraphics[width=0.48\textwidth]{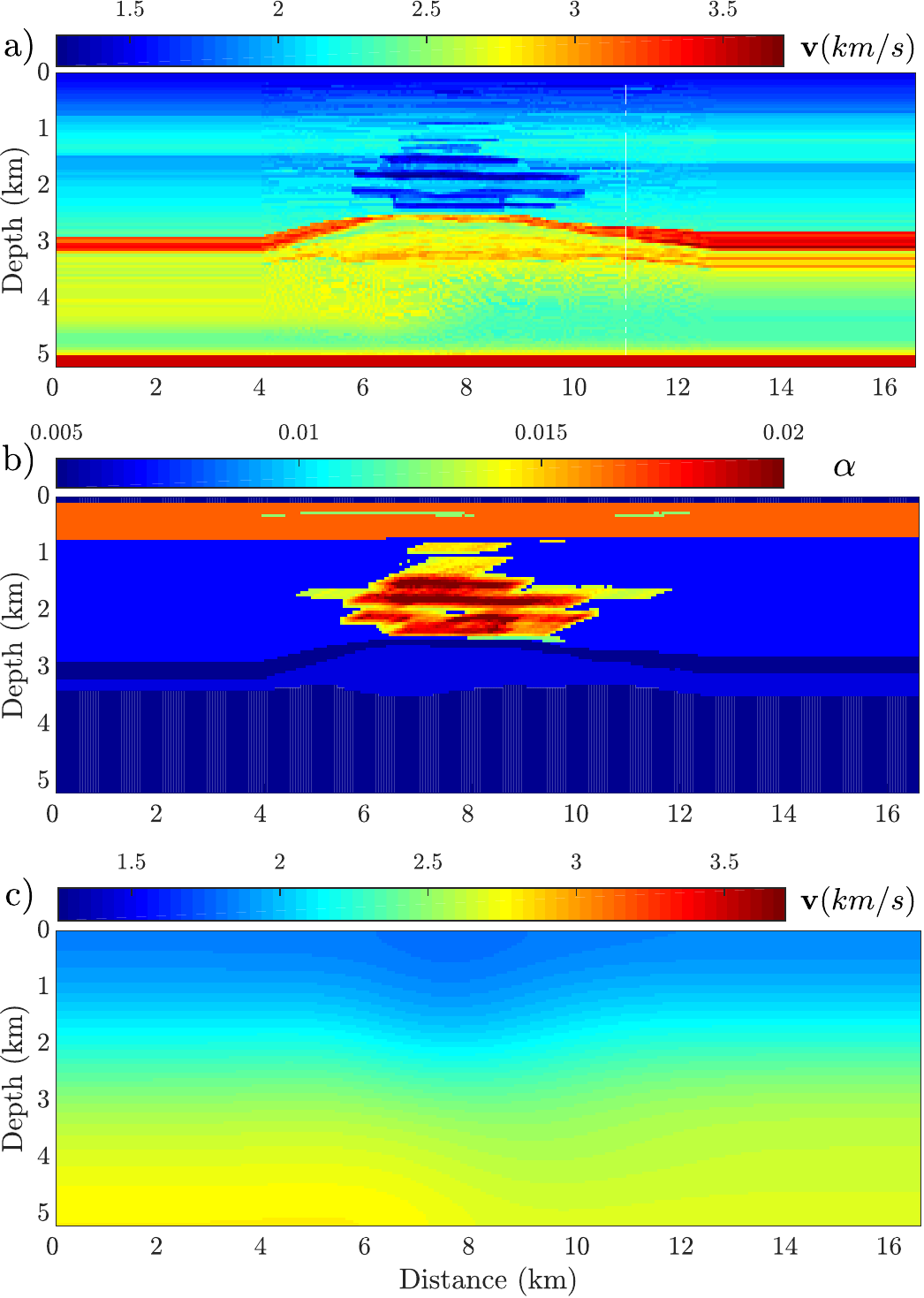}
\caption{North Sea case study. (a) True $\bold{v}$ model. (b) True $\boldsymbol{\alpha}$ model. (c) Initial $\bold{v}$ model.}
\label{fig:north_test_true}
\end{figure}

We consider a more realistic 16~km $\times$ 5.2~km shallow-water synthetic model representative of the North Sea \citep{Munns_1985_VFG}. The true $\bold{v}$ and $\boldsymbol{\alpha}$ models are shown in Figures \ref{fig:north_test_true}a and \ref{fig:north_test_true}b, respectively. 
The velocity model is formed by soft sediments in the upper part, a pile of low-velocity gas layers above a chalk reservoir, the top of which is indicated by a sharp positive velocity contrast at around 2.5~km depth, and a flat reflector at 5~km depth (Figure \ref{fig:north_test_true}a). 
The $\boldsymbol{\alpha}$ model has two highly attenuative zones in the upper soft sediments and gas layers, and the $\boldsymbol{\alpha}$ value is relatively low elsewhere (Figure \ref{fig:north_test_true}b). 

The fixed-spread surface acquisition consists of 320 explosive sources spaced 50~m apart at 25~m depth and 80 hydrophone receivers spaced 200~m apart on the sea floor at 75~m depth. For sake of computational efficiency, we use the spatial reciprocity of Green's functions to process sources as receivers and vice versa.
A free-surface boundary condition is used on top of the grid and the source signature is a Ricker wavelet with a 10~Hz dominant frequency.
We perform forward modeling with a nine-point stencil finite-difference method implemented with anti-lumped mass and PML absorbing boundary conditions, where the stencil coefficients are optimized to the frequency \citep{Chen_2013_OFD}. We solve the normal-equation system for wavefield reconstruction, equation~\ref{eqnormal}, with a sparse direct solver \citep{Duff_2017_DMS}.
%
%
%
%
%
\begin{figure*}
\includegraphics[width=1\textwidth]{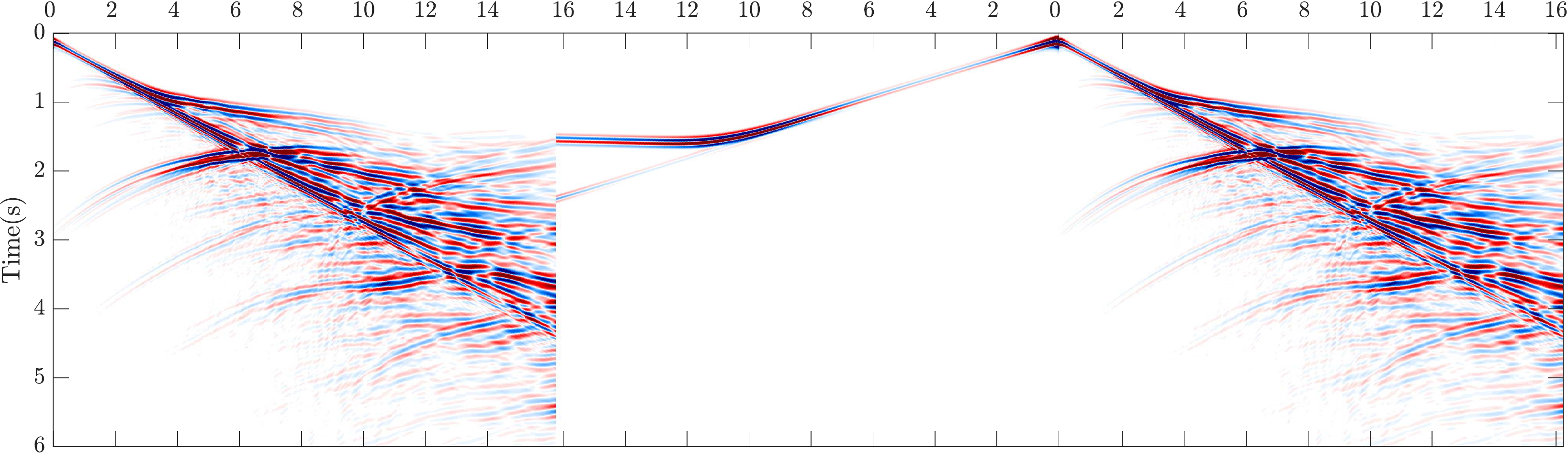}\\
\caption{Time domain seismograms computed in the true and initial models. The true seismograms are shown in the left and right panels, while those computed in the initial model are shown in the middle panel with a mirror representation such that the two sets of seismograms can be compared at long and short offsets.  The seismograms are plotted with a reduction velocity of 2.5 km/s for sake of time axis compression.}
\label{fig:north_test_seis_t}
\end{figure*}
%
%
The initial model for $\bold{v}$ is a highly Gaussian filtered version of the true model (Figure \ref{fig:north_test_true}c), while the starting $\boldsymbol{\alpha}$ model is homogeneous with $\boldsymbol{\alpha}=0$. The common-shot gathers computed in the true and initial models are compared in Figure \ref{fig:north_test_seis_t} for a shot located at 16.0~km. 
The latter mainly show the direct wave and the diving waves, which are highly cycle skipped relative to those computed in the true model. 

We perform the inversion with small batches of three frequencies with one frequency overlap between two consecutive batches, moving from the low frequencies to the higher ones according to a classical frequency continuation strategy. The starting and final frequencies are 3~Hz and 15~Hz and the sampling interval in each batch is 0.5~Hz. We perform three paths through the frequency batches to improve the results, using the final model of one path as the initial model of the next one (these paths can be viewed as outer iterations of the algorithm). The starting and finishing frequencies of the paths are [3, 6], [4, 10], [6, 15]~Hz respectively, where the first and second elements of each pair show the starting and finishing frequencies, respectively. Also, we used batches of four frequencies with two frequencies overlap during the second path and five frequencies with three frequencies overlap during the third path.
The motivation behind this frequency management is to keep the bandwidth of each patch narrow during the first path in order to mitigate nonlinearities through a progressive frequency continuation, before broadening this bandwidth during the second and third paths to strengthen the imprint of dispersion in the inversion and decouple velocity and attenuation more efficiently.

\subsubsection{Comparison of FWI and WRI objective Functions}
\begin{figure}[ht!] 
\begin{center}
  \includegraphics[width=0.48\textwidth]{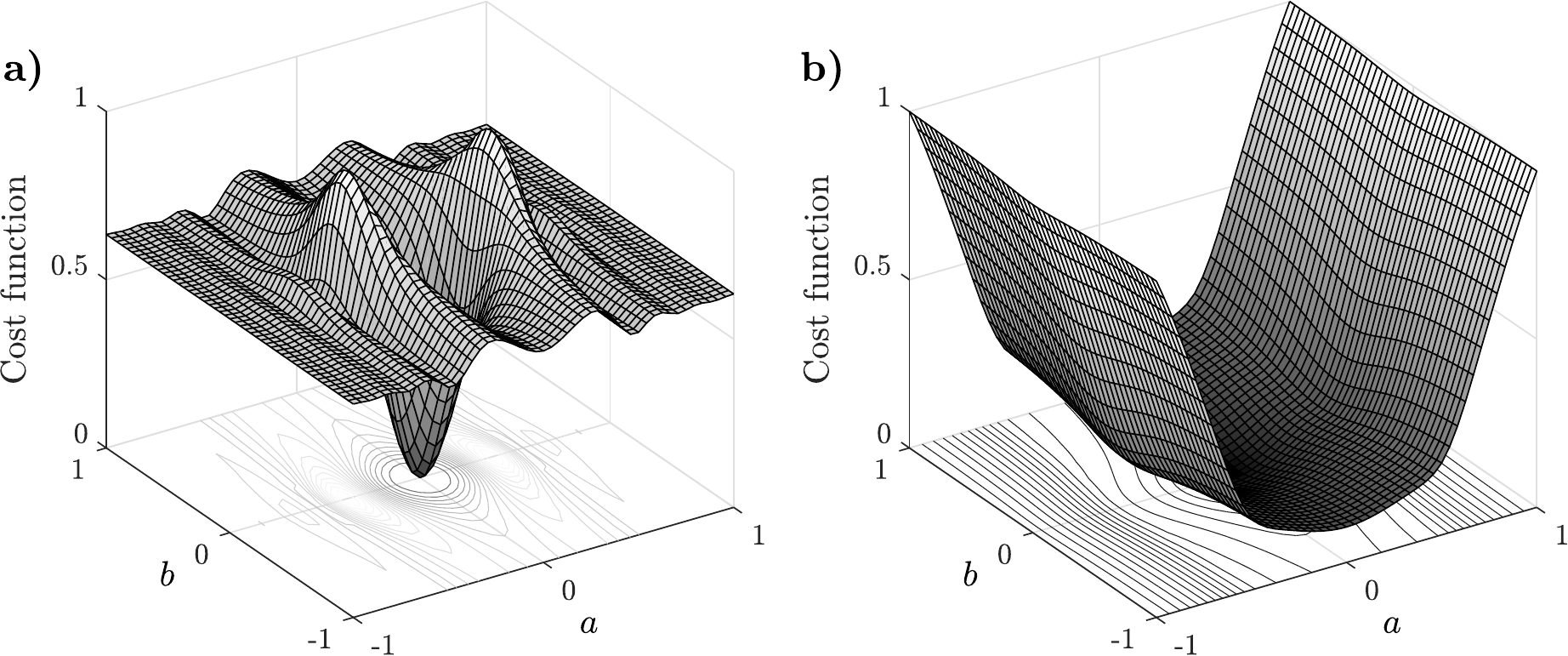} 
\caption{The objective function for the 3 Hz wavefield as a function of $\bold{v}_{a}$ and $\boldsymbol{\alpha}_b$ generated using equation \ref{eq_cost}. (a) Classical reduced-space FWI. (b) WRI.}
\label{fig:val_cost}
\end{center}
\end{figure}

Before showing the inversion results, it is worth illustrating how WRI extends the search space of FWI for this North Sea case study. For this purpose,  we compare the shape of the classical FWI misfit function based upon the $\ell{2}$ norm of the data residuals \citep[e.g. ][]{Pratt_1998_GNF} with that of the parameter-estimation WRI subproblem for the 3~Hz frequency and for a series of $\bold{v}$ and $\boldsymbol{\alpha}$ models that are generated according to 
\begin{subequations} \label{eq_cost}
\begin{eqnarray}
\bold{v}_{a}=\bold{v}_{true} + a^2 (\bold{v}_{init}-\bold{v}_{true}), \\
\boldsymbol{\alpha}_b=\boldsymbol{\alpha}_{true} + b^2 (\boldsymbol{\alpha}_{init}-\boldsymbol{\alpha}_{true}),
\end{eqnarray}
\end{subequations}
where $-1 \leq a,b \leq 1$. In this case, $(\bold{v}_{a},\boldsymbol{\alpha}_b)$ lies on the line-segment joining the initial point ($\bold{v}_{init},\boldsymbol{\alpha}_{init})$ and final point ($\bold{v}_{true},\boldsymbol{\alpha}_{true})$.  
We set $\bold{v}_{true}$ and $\boldsymbol{\alpha}_{true}$ as those shown in Figures \ref{fig:north_test_true}(a-b), 
 $\bold{v}_{init}$ as the initial $\bold{v}$ model shown in Figure \ref{fig:north_test_true}c, and $\boldsymbol{\alpha}_{init}$ as a homogeneous model with $\boldsymbol{\alpha} = 0.004$ ($\bold{Q}=250$).
The misfit functions of FWI and WRI are shown in Figures \ref{fig:val_cost}a and \ref{fig:val_cost}b, respectively. 
The FWI objective function exhibits spurious local minima with respect to both velocity and attenuation ($a$ and $b$ dimensions), while only one minimum is seen in the WRI objective function. Indeed, this highlights the search space expansion generated by the wave-equation relaxation during WRI. This search space expansion is displayed through a wider and flatter attraction basin compared to that of FWI. This wide attraction basin exacerbates the contrast between the sensitivities of the objective function to velocity and attenuation, with a quite weak sensitivity to the latter. These contrasted sensitivities would likely make the parameter estimation subproblem poorly scaled if the velocity and attenuation were jointly updated during WRI based on variable projection \citep{vanLeeuwen_2016_PMP}. In this context, the alternating-direction strategy can be viewed as an heuristic to overcome this scaling issue. Indeed, the lack of sensitivity to attenuation during the early WRI iterations requires aggressive regularization and bound constraints to stabilize the attenuation estimation. As WRI proceeds over iterations and the wave equation constraint is satisfied more accurately, the inversion should recover a significant sensitivity to attenuation as that highlighted in Figure~\ref{fig:val_cost}a allowing for a relaxation of the regularization. These statements highlight the need to reconcile search space expansion to manage nonlinearity and regularization plus bound constraints to restrict the range of feasible solutions for $\boldsymbol{\alpha}$. \\

\subsubsection{Viscoacoustic FWI results}
\begin{figure}[htb!]
\begin{center}
\includegraphics[width=0.48\textwidth]{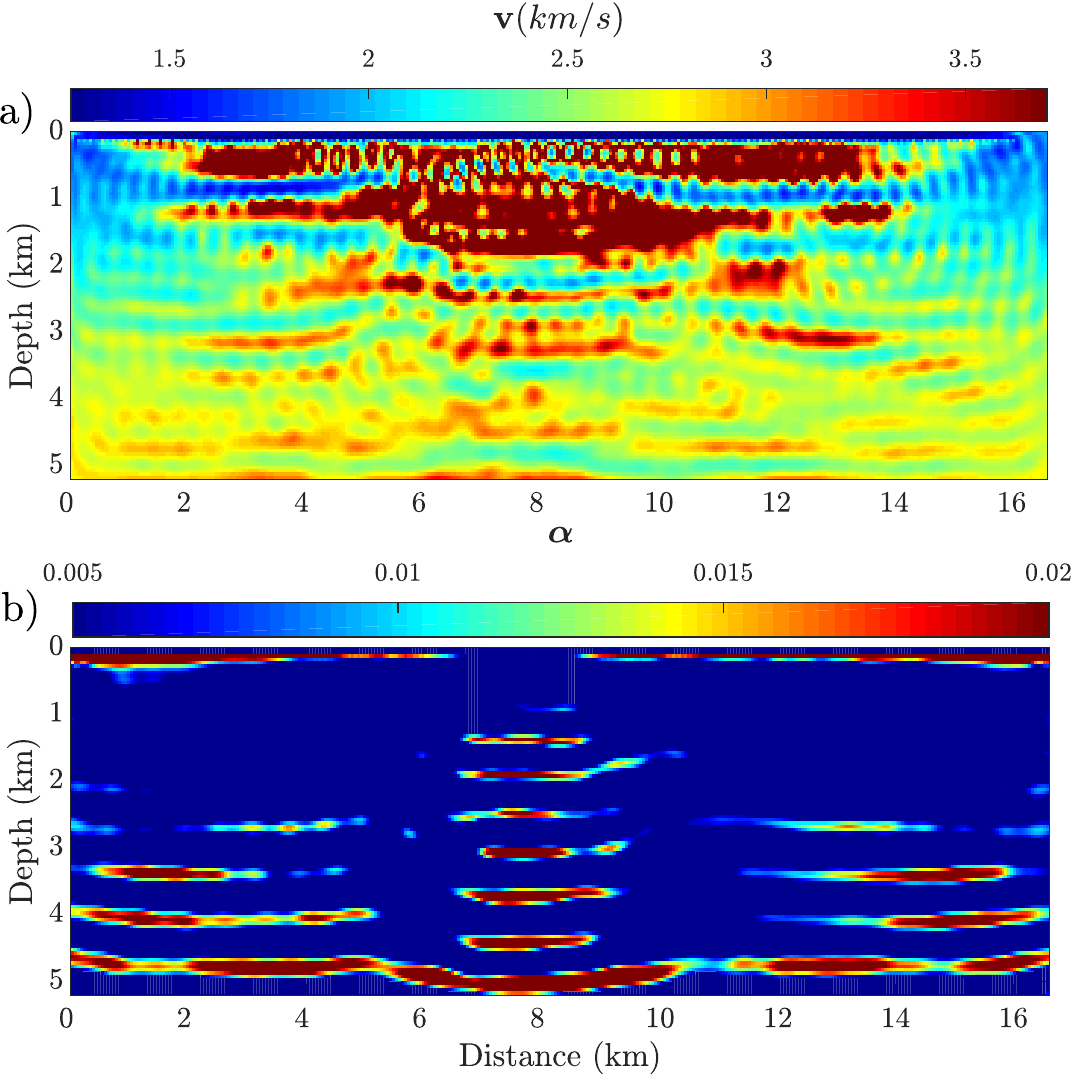}\\
\caption{Viscoacoustic FWI results after inverting the first frequency batch. (a) Reconstructed velocity model. (b) Reconstructed attenuation model. TV regularization and bound constraints are not applied.}
\label{fig:north_test_red}
\end{center}
\end{figure}

To show the need of search space expansion and regularization, we first perform a classical viscoacoustic FWI for noiseless data \citep[e.g. ][]{Kamei_2013_ISV,Operto_2018_MFF}. We use squared-slowness and attenuation as optimization parameters and update them simultaneously with the L-BFGS quasi-Newton optimization and a line search procedure for step length estimation (that satisfies the Wolfe conditions). Neither TV regularization nor bound constraints are applied. 
Owing the limited kinematic accuracy of the initial models highlighted by the seismograms mismatches in Figure \ref{fig:north_test_seis_t}, the reconstruction of the velocity model remains stuck in a local minimum during the first frequency batch inversion (Figure~\ref{fig:north_test_red}a), while the estimated attenuation model shows unrealistic values due to the lack of bound constraints (Figure~\ref{fig:north_test_red}b). This failure prompts us to stop the inversion at this stage.

\subsubsection{Viscoacoustic IR-WRI results}
\begin{figure*}[ht!]
\includegraphics[width=1\textwidth]{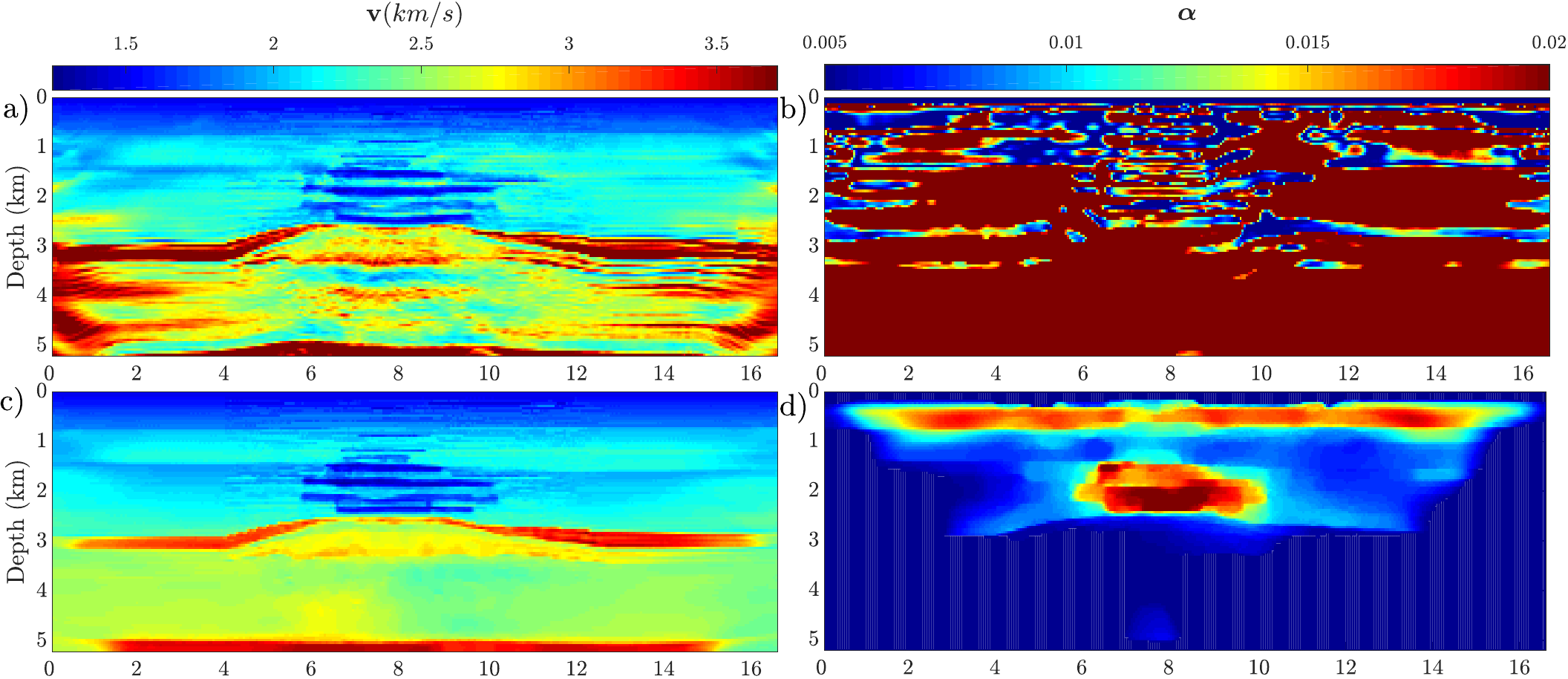}\\
\includegraphics[width=0.99\textwidth]{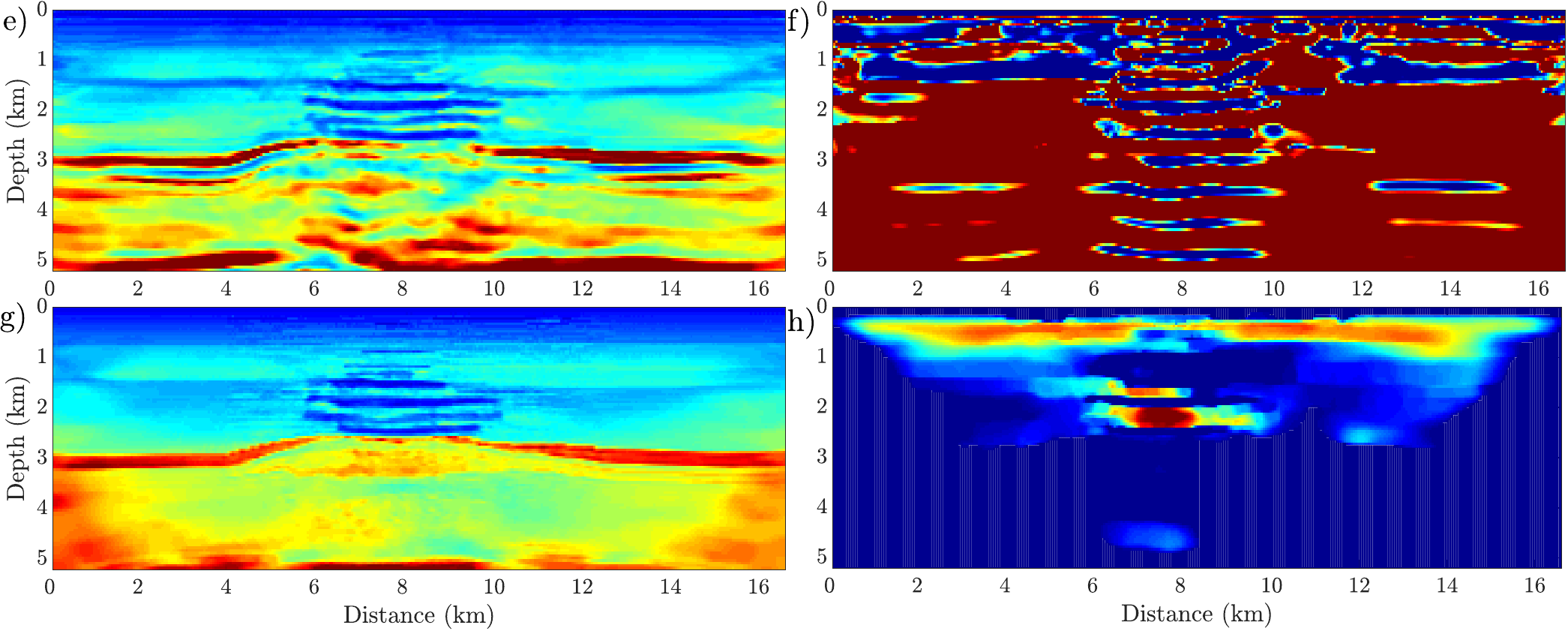}\\
\caption{Viscoacoustic IR-WRI results. (a-d) Noiseless data. (a-b) Velocity (a) and attenuation (b) models reconstructed without BTV regularization. (c-d) Same as (a-b) when BTV regularization is applied. (e-h) Same as (a-d) for noisy data  (SNR=10db).}
\label{fig:north_test}
\end{figure*}

\begin{figure}[ht!]
\centering
\includegraphics[width=0.49\textwidth]{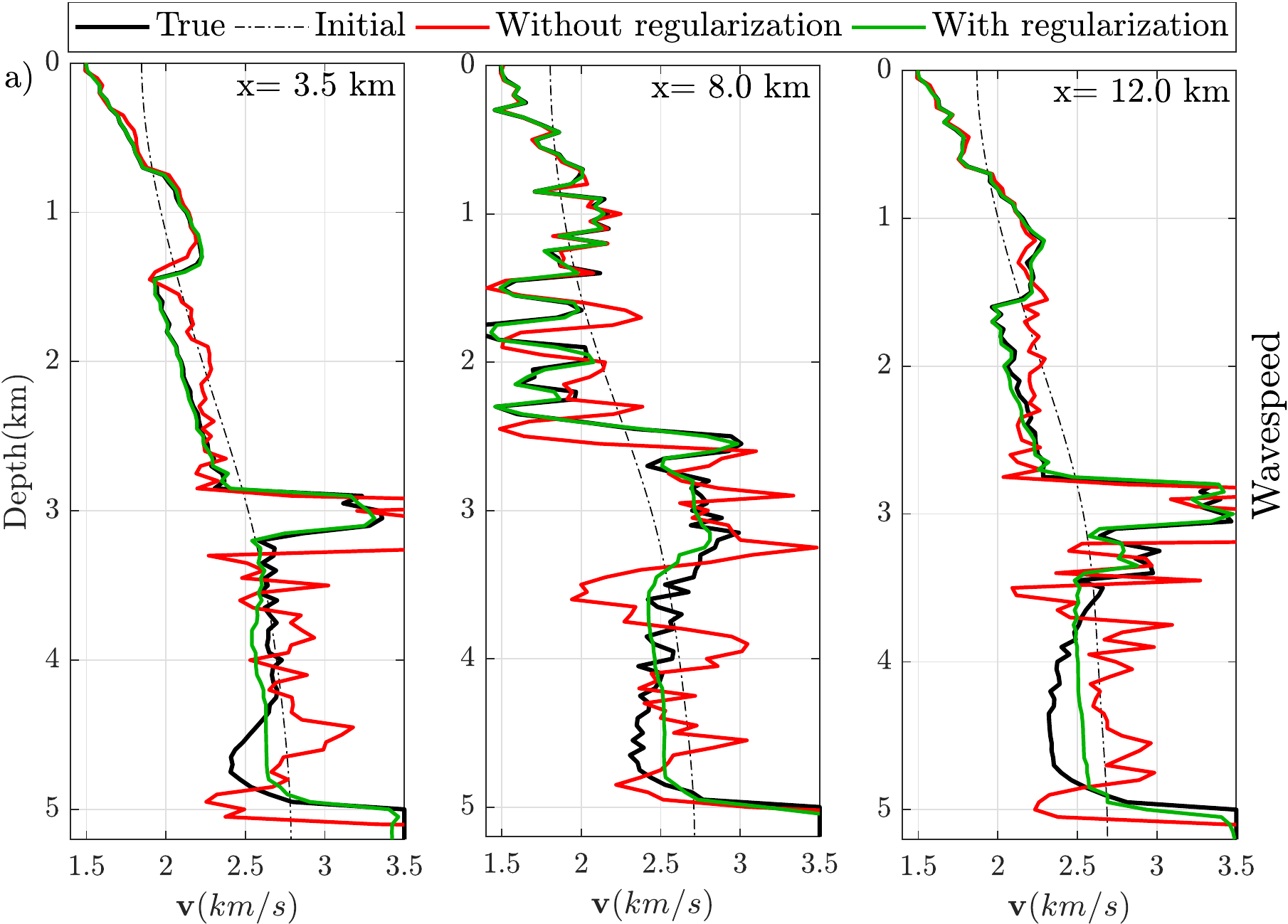}\\
\includegraphics[width=0.49\textwidth]{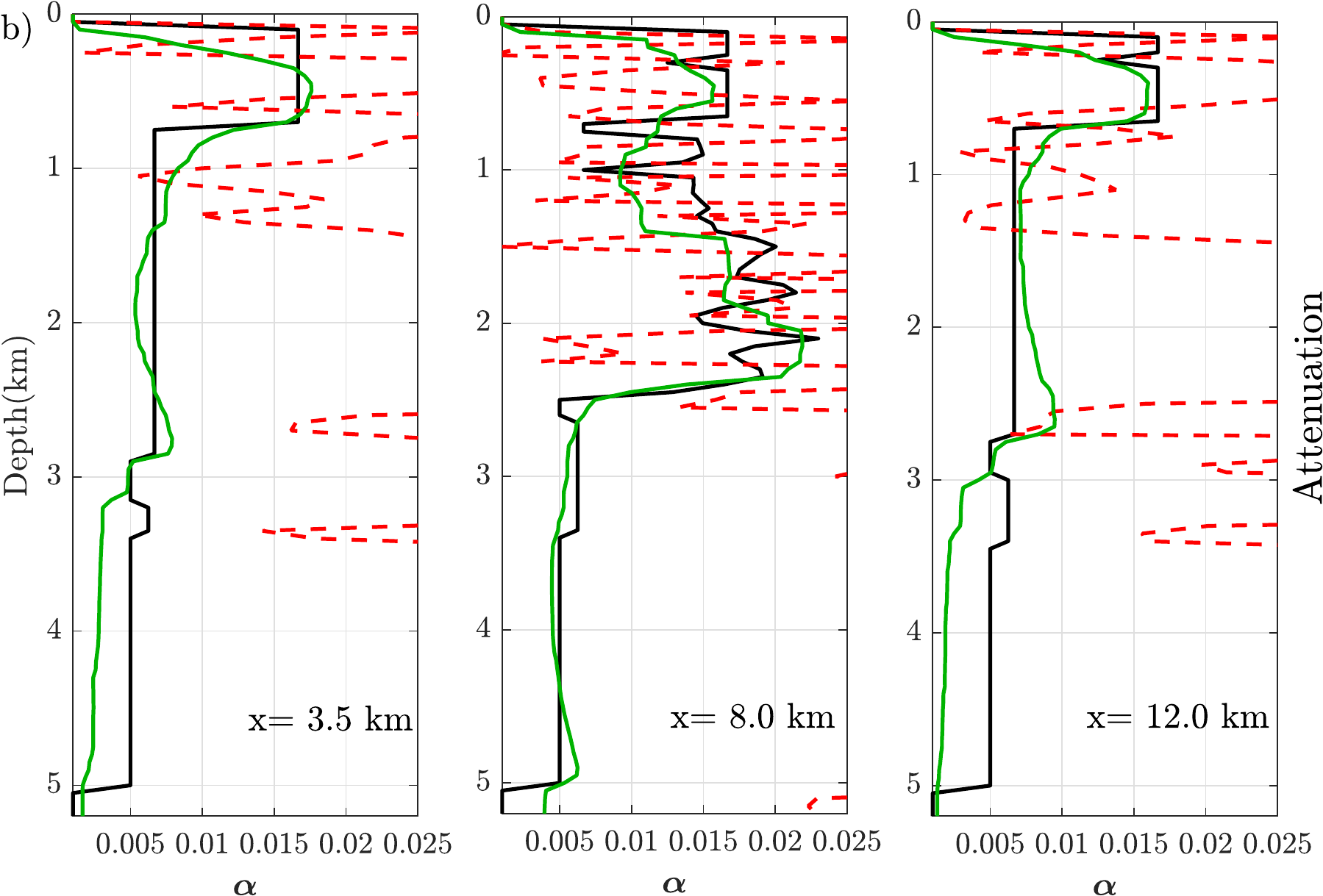}
\caption{For noiseless data, direct comparison along the logs at $x=3.5$ (left), $x=8.0$ (center) and $x=12.0$ km (right) between the true model (black), the initial model (dashed line) and the IR-WRI models without (red) and with BTV  (green) regularization (Figure~\ref{fig:north_test}(a-d)). (a) Estimated $\bold{v}$ , (b) estimated $\boldsymbol{\alpha}$.}
\label{fig:north_test_log}
\end{figure}

We now update $\bold{v}$ and $\boldsymbol{\alpha}$ with IR-WRI according to the optimization workflow described in Algorithm 1. IR-WRI is performed without and with bound constraints + TV regularization (referred to as BTV regularization in the following). The lower and upper bounds are 1.2~km/s and 4~km/s for velocities, and 0.001 and 0.025 for $\boldsymbol{\alpha}$.
For each case, the stopping criterion for each batch is given to be either reaching a maximum iteration count of 20 or 
\begin{align}
\label{Stop}
&\sum\| \bold{A}(\bold{m}^{k+1},\boldsymbol{\alpha}^{k+1}) \bold{u}^{k+1}-\bold{b}\|_2^2 \leq \epsilon_b ~~ \nonumber ~~~\text{and}\\
 & \sum\|\bold{Pu}^{k+1}-\bold{d}\|_2^2 \leq \epsilon_d,
\end{align} 
where the sums run over the frequencies of the current batch, $\epsilon_b$=1e-3 and $\epsilon_d$=1e-5.
%
%
%
We start with inversion of noiseless data. The final $\bold{v}$ and $\boldsymbol{\alpha}$ models, estimated by IR-WRI without and with BTV regularization, are shown in Figure \ref{fig:north_test}(a-d) after 360 and 321 iterations, respectively. A direct comparisons between the logs extracted from the true models, the initial model and the IR-WRI velocity models reconstructed without/with BTV regularization at $x=3.5$~km, $x=8.0$~km and $x=12.0$~km are shown in Figure \ref{fig:north_test_log}a. 
Although a crude initial velocity model was used, the velocities in the shallow sedimentary cover and the gas layers are fairly well reconstructed in both cases (Figure \ref{fig:north_test}a,c). Also, IR-WRI without BTV regularization manages to reconstruct an acceptable velocity model unlike classical FWI (Figure~\ref{fig:north_test_red}).  The main differences between the IR-WRI velocity models built with and without BTV regularization are shown at the reservoir level and below. Without BTV regularization, the top of the reservoir is mispositioned (Figure~\ref{fig:north_test_log}a, $x=8.0$~km, green versus red curves) and the inversion fails to reconstruct the smoothly-decreasing velocity below it due to the lack of diving wave illumination at these depths (Figure~\ref{fig:north_test_log}a). This in turn prevents the focusing of the deep reflector at 5~km depth by migration of the associated short-spread reflections (Figure~\ref{fig:north_test}a). When BTV regularization is used, viscoacoustic IR-WRI provides a more accurate and cleaner images of the reservoir and better reconstructs the sharp contrast on top of it (Figure~\ref{fig:north_test}c). As expected, the TV regularization replaces the smoothly-varying velocities below the reservoir (between 3 to 5 km depth) by a piecewise-constant layer due to the lack of wave illumination in this part of the model (Figure~\ref{fig:north_test_log}a, green curves). However, this does not prevent a fairly accurate reconstruction of the deep reflector at 5~km depth.

A direct comparisons between the logs extracted from the true and the IR-WRI attenuation models at $x=3.5$~km, $x=8.0$~km and $x=12.0$~km are shown in Figure \ref{fig:north_test_log}b. 
The reconstruction of $\boldsymbol{\alpha}$ without BTV regularization is quite unstable, with an oscillating trend and overestimated values  (Figures \ref{fig:north_test}b and red curves in \ref{fig:north_test_log}b). This highlights fairly well the ill-posedness of the attenuation reconstruction.
In contrast, the $\boldsymbol{\alpha}$ model reconstructed with BTV regularization captures the large-scale attenuation trend in the shallow sedimentary cover, in the gas layers and below the reservoir (Figures \ref{fig:north_test}d and green curves in \ref{fig:north_test_log}b). We note however that the attenuation is underestimated on top of the gas layers between 1~km and 1.4~km depth (Figure~\ref{fig:north_test_log}b, $x=8.0$~km, green curve). This error in the attenuation reconstruction might be correlated with subtle underestimation of velocities at these depths (Figure~\ref{fig:north_test_log}a, $x=8$~km, green curve). This might indicate on the one hand some mild amplitude-related cross-talk effects between velocities and attenuation and on the other hand the higher sensitivity of the data to velocities compared to attenuation (in the sense that a small error in the velocity contrasts can compensate more significant attenuation errors). Similar cross-talk artifacts have been previously discussed during the toy inclusion test (Figure~\ref{fig:Box_test}e-f).

\begin{figure}[ht!]
\centering
\includegraphics[width=0.49\textwidth]{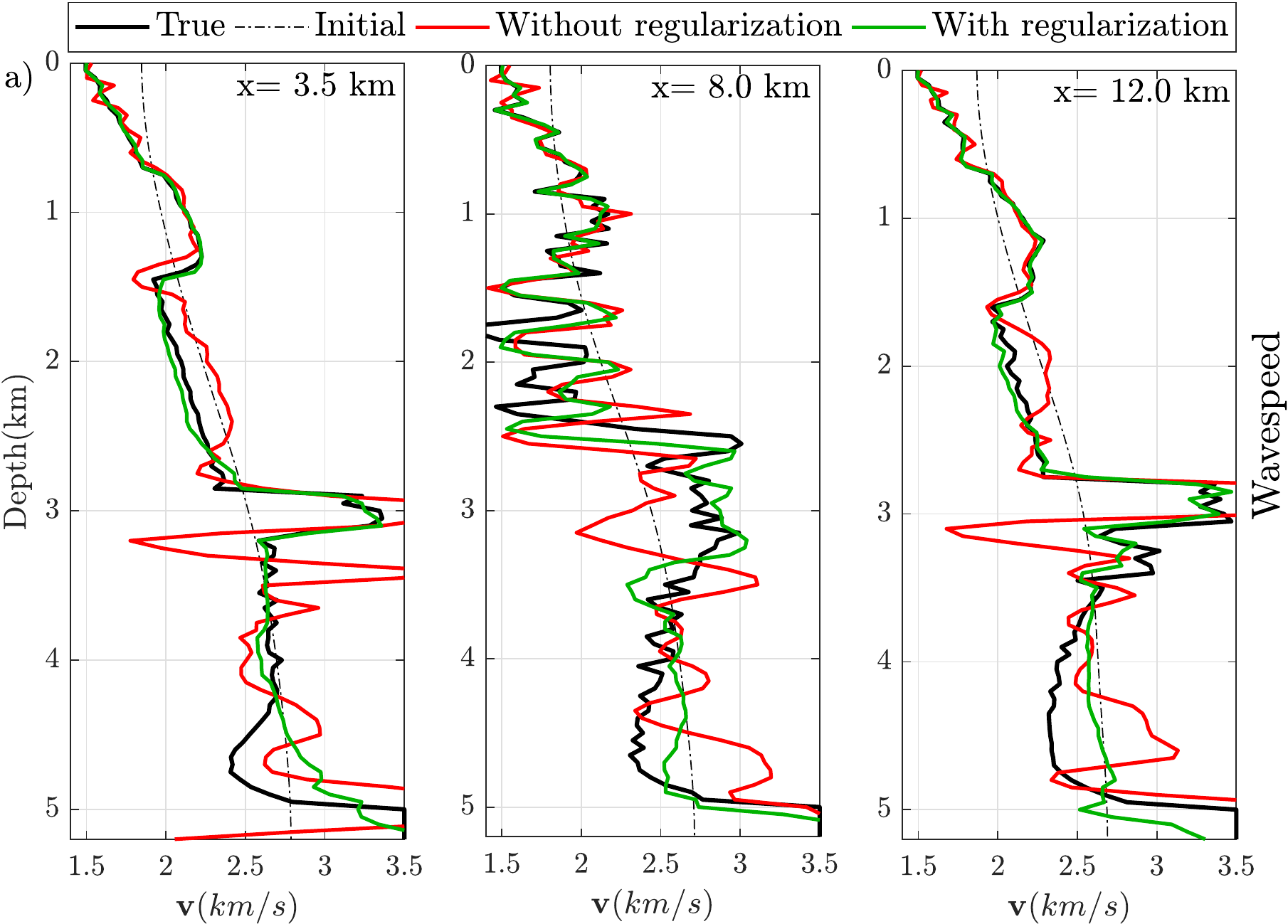}\\
\includegraphics[width=0.49\textwidth]{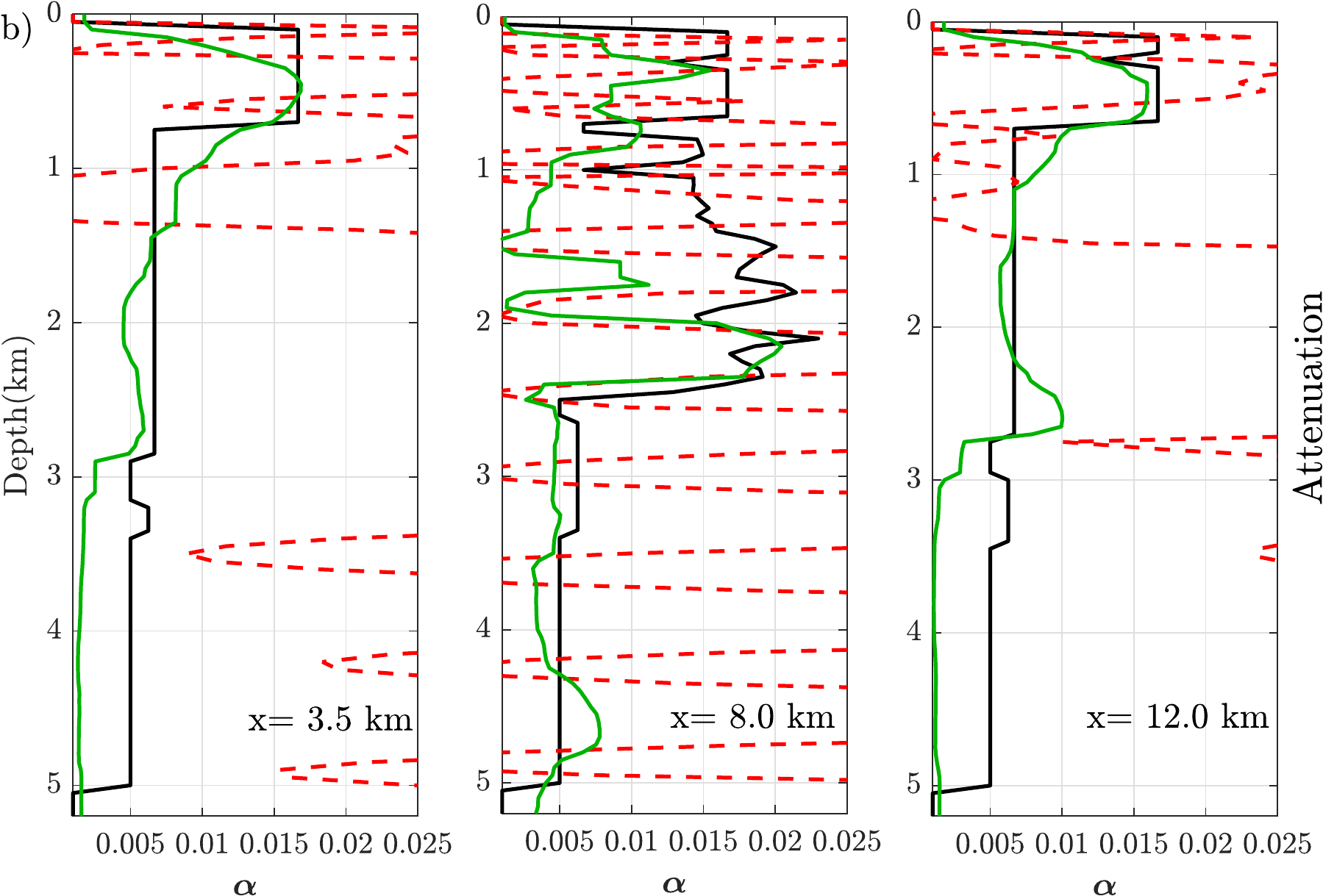}
\caption{Same as Fig. \ref{fig:north_test_log} but for noisy data with SNR=10db. The corresponding IR-WRI models are shown in Figure~\ref{fig:north_test}(e-h)).}
\label{fig:north_test_log1}
\end{figure}

\begin{figure*}[ht!]
\includegraphics[width=1\textwidth]{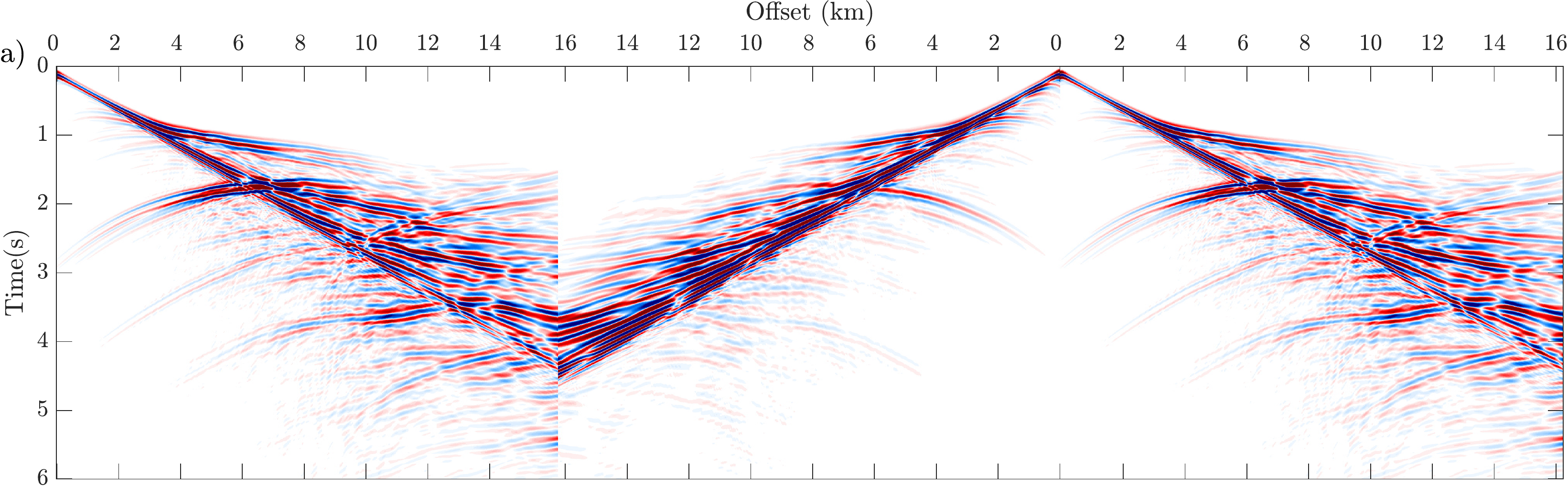}\\
\includegraphics[width=1\textwidth]{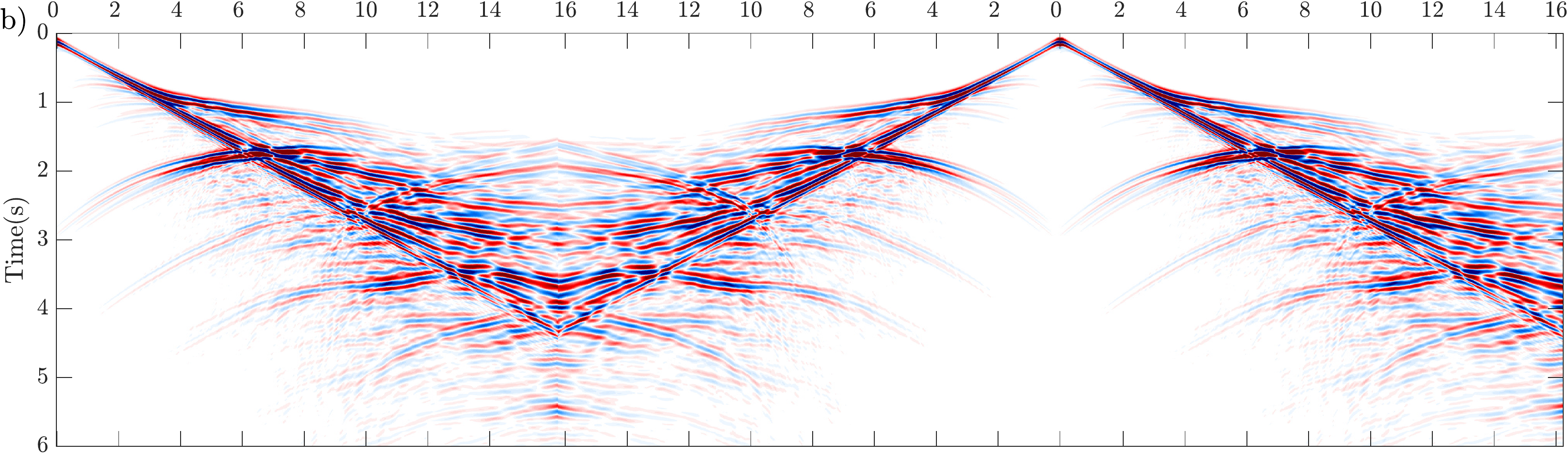}\\
\includegraphics[width=1\textwidth]{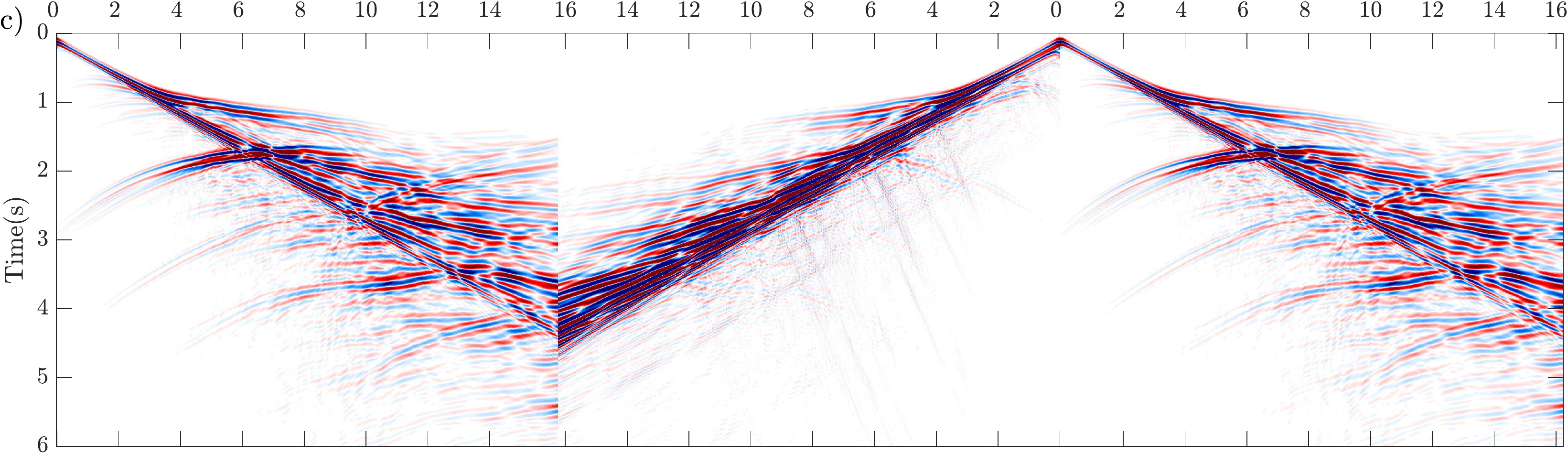}\\
\includegraphics[width=1\textwidth]{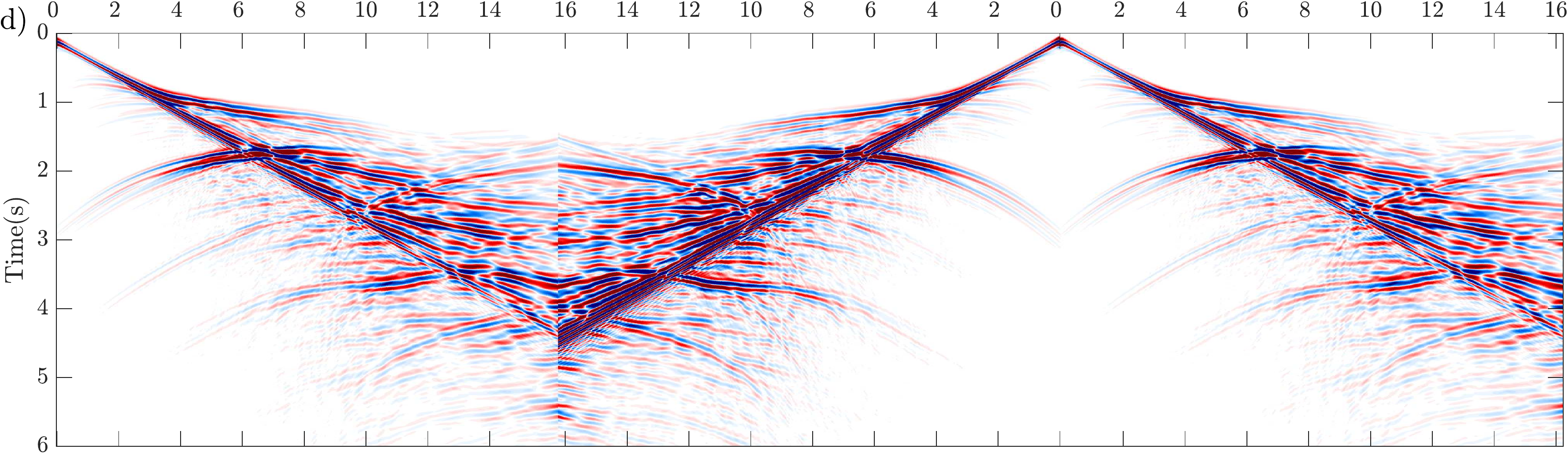}
\caption{Time domain seismograms computed in (a) IR-WRI without regularization and (b) BTV regularized IR-WRI for noiseless data (Figure \ref{fig:north_test}a- \ref{fig:north_test}d). (c-d) Same as (a-b), but for noisy data (Figure \ref{fig:north_test}e- \ref{fig:north_test}h).  The true seismograms are shown in the first and the last panel of the above mentioned seismograms (folded) to have a comparison at short and long offset with true seismograms.  The seismograms are plotted with a reduction velocity of 2.5 km/s for sake of time axis compression.}
\label{fig:north_test_seis}
\end{figure*}

We continue by assessing the resilience of the proposed viscoacoustic IR-WRI to noise when data are contaminated with a Gaussian random noise with a SNR=10~db. Here, SNR is defined based on the root mean square (RMS) amplitude of signal and that of noise as 
\begin{equation}
\text{SNR}=20 \log \left( \frac{Signal ~RMS~ Amplitude}{Noise ~RMS ~Amplitude} \right).
\end{equation}
We use the same setup and the same initial models as those used for the noiseless case. The stopping criterion is defined by equation \ref{Stop}, where $\varepsilon_d$ is now set to the noise level. 
The final models of IR-WRI obtained without and with BTV regularization are shown in Figure \ref{fig:north_test}(e-h). The total number of IR-WRI iterations are 196 and 185, respectively, for these results. In a similar manner to the noiseless case,  a direct comparisons between the logs extracted from the true models, the initial model and the IR-WRI velocity models at $x=3.5$~km, $x=8.0$~km and $x=12.0$~km are shown in Figure \ref{fig:north_test_log1}. Overall, a similar trend as for the noiseless case is shown. However, the presence of noise in the data leads to a mispositioning of the reservoir at 8~km distance in the BTV IR-WRI velocity model, which was not observed in the noiseless case (compare Figures~\ref{fig:north_test_log}a and \ref{fig:north_test_log1}a, x= 8km, green curves). This mispositioning of the reservoir may be correlated with a poorer reconstruction of the attenuating gas layers between 1~km and 2.5~km depth (compare Figures~\ref{fig:north_test_log}b and \ref{fig:north_test_log1}b, x= 8km, green curves). \\
We also compute a common-shot gather for a shot located at 16.0~km in the IR-WRI models inferred from noiseless/noisy data with/without BTV regularization (Figure \ref{fig:north_test_seis}).  
The time-domain seismograms computed in the IR-WRI models obtained without regularization (Figure \ref{fig:north_test_seis}a and \ref{fig:north_test_seis}c) show underestimated amplitudes and do not match late dispersive arrivals due to the overestimated and oscillating values of $\boldsymbol{\alpha}$ (Figures \ref{fig:north_test}b and \ref{fig:north_test_log}b for noiseless data and Figures \ref{fig:north_test}f and \ref{fig:north_test_log1}b for noisy data). \\
In the case of noiseless data, the bound constraints and the TV regularization allow for a high-quality data fit (Figure \ref{fig:north_test_seis}c), consistently with the accuracy of the models shown in Figure~\ref{fig:north_test}c-d. In the case of noisy data, the bound constraints and the TV regularization improve significantly the data match (compare Figure \ref{fig:north_test_seis}c and Figure \ref{fig:north_test_seis}d). However, the imprint of the cross-talk artifacts mentioned above are clearly seen at long offsets with a degraded fit of deeply-propagating waves (for example, the refracted wave from the deep reflector at around 1.5s~traveltime and the late dispersive waves at around 4~s traveltime) relative to the noiseless-data results (compare Figure \ref{fig:north_test_seis}b and Figure \ref{fig:north_test_seis}d).


%
%
%
%

\section{Conclusions} 
We extended the recently proposed ADMM-based iteratively-refined wavefield reconstruction inversion (IR-WRI) for attenuation imaging by inversion of viscoacoustic wavefields.
The proposed viscoacoustic IR-WRI treats the nonlinear viscoacoustic waveform inversion as a multiconvex optimization problem. To achieve this goal, the original nonlinear multi-parameter problem for squared slowness and attenuation factor is replaced by three recursive linear mono-parameter subproblems for wavefield, squared slowness and attenuation factor that are solved in alternating mode at each IR-WRI iteration. The attenuation-reconstruction subproblem requires to introduce an approximate multilinear viscoacoustic wave equation in wavefield, squared slowness, and attenuation factor. However, the errors generated by this approximate viscoacoustic wave equation during the attenuation reconstruction are efficiently compensated by the Lagrange multipliers (namely, the running sum of the wave equation errors) that are computed with the exact viscoacoustic equation. \\
This new formulation has first the flexibility to tailor the regularization to the squared slowness and attenuation factor. Moreover, it simplifies the multi-parameter optimization workflow and mitigates its computational cost since the original poorly-scaled multi-parameter inversion is recast as two interlaced mono-parameter inversions. 
A realistic synthetic example suggests that the search space extension embedded in IR-WRI efficiently mitigates cycle skipping when a crude initial velocity model is used, while the TV-regularized alternating-direction optimization reasonably manages the cross-talks between squared slowness and attenuation as well as the limited sensitivity of the data to the attenuation.
%
%


\append{Scaled form of augmented Lagrangian}
In this appendix, we briefly review how augmented Lagrangian (AL) function, as the one shown in equation~\ref{eqpsi}, is used to solve constrained problem with the method of multiplier \citep[][ Chapter 17]{Nocedal_2006_NOO}.
Let's start with the following constrained problem
\begin{eqnarray}\label{app1}
\min_{\bold{x}} ~ \|P(\bold{x})\|_2^2  ~~~~ \text{subject to} ~~~~Q(\bold{x})=\bold{0}.
\end{eqnarray} 
The AL function associated with the problem \ref{app1} combines a Lagrangian function and a penalty function as
\begin{equation} \label{app2}
\mathcal{L}_A(\bold{x},\bold{v}) = \underbrace{\|P(\bold{x})\|_2^2 +  \langle\bold{v}, Q(\bold{x})\rangle}_{\text{Lagrangian}} + 
\underbrace{\frac{\xi}{2}\|Q(\bold{x})\|_2^2}_{\text{Augmentation}},
\end{equation}
where $\langle \cdot,\cdot\rangle$ denotes inner product and $\bold{v}$ and $\xi$ denotes the Lagrange multiplier (dual variable) and the penalty parameter, respectively. \\
This AL function can be written in a compact form by introducing the scaled dual variable ${\bold{q}}=-\bold{v}/\xi$ and adding and subtracting the term $\frac{\xi}{2}\|{\bold{q}}\|_2^2$ to the function \ref{app2} \citep[See ][ Page 15 for more details]{Boyd_2011_DOS}:
\begin{align}\label{app3}
\mathcal{L}_A(\bold{x},{\bold{q}}) &=\|P(\bold{x})\|_2^2 - \xi \langle\bold{q}, Q(\bold{x})\rangle+ \frac{\xi}{2}\|Q(\bold{x})\|_2^2+\frac{\xi}{2}\|{\bold{q}}\|_2^2-\frac{\xi}{2}\|{\bold{q}}\|_2^2 \nonumber \\
&=\|P(\bold{x})\|_2^2+ \frac{\xi}{2} \|Q(\bold{x})-{\bold{q}}\|_2^2 - \frac{\xi}{2}\|{\bold{q}}\|_2^2.
\end{align} 
Equation~\ref{app3} shows the augmented Lagrangian method can be seen as a penalty method with an error correction term in the penalty function, $\frac{\xi}{2} \|Q(\bold{x})-{\bold{q}}\|_2^2$, corresponding to the scaled Lagrange multipliers. This correction term controls how well the constraint is satisfied at the convergence point.
In the framework of the method of multiplier, the AL function is minimized with respect to the primal variable $\bold{x}$ and maximized with respect to the scaled dual variable ${\bold{q}}$ in alternating mode.
Expression \ref{app3} shows that the dual variable is simply updated with the constraint violation when a gradient ascent method is used. This recipe has been used to derive equation~\ref{ADMM} with ADMM (AL method with alternating update of multiple classes of primal variable). The reader is referred to \citet{Aghamiry_2019_IWR} for the detailed development.
%
\append{Bound constrained  TV-regularization using ADMM} \label{Appa}
\label{appendB}
In this appendix, we review step by step how to solve a bound-constrained TV-regularized convex problem (such as those in equations \ref{TV_sig} and \ref{TV_eta}) using variable splitting and ADMM \citep{Boyd_2011_DOS}. 
Let's consider a general bound-constrained TV-regularized convex problem of the form 
\begin{equation}  \label{CTV}
 \underset{\bold{x}\in \mathcal{X}}{\min}~
 \sum_{i=1}^N  \sqrt{|\nabla_{\!x} \bold{x}|_i^2 + |\nabla_{\!z}\bold{x}|_i^2} + \frac{\lambda}{2} \|\bold{G}\bold{x}-\bold{y}\|_2^2,
\end{equation}
for some column vector $\bold{y}$ and matrix $\bold{G}$. The model $\bold{x}$ is an $N$-length column vector, $\nabla_{\!x}$ and $\nabla_{\!z}$ are square first-order difference matrices, and $\mathcal{X}$ is the desired convex set. The penalty parameter $\lambda>0$ balances the relative weight of the TV regularizer and the misfit term.
Following \citet[ section 2.2.2]{Aghamiry_2019_IBC}, equation \ref{CTV} can be solved via the following three easy tricks.\\ 

1) \textit{Variable splitting.} Since the variable $\bold{x}$ appears simultaneously in the TV, misfit, and bounding terms, it ``couples" these terms and makes it difficult to solve the problem. To decouple them, new auxiliary variables $\bold{p}_x=\nabla_{\!x} \bold{x}$, ${\bold{p}_y}\in \mathcal{X}$, and $\bold{p}_z=\nabla_{\!z} \bold{x}$ are substituted in the TV term and the bound constraint, respectively, and their expression as a function of the original variable $\bold{x}$ are introduced as new equality constraints. This recasts \ref{CTV} as the following constrained problem
\begin{align}  \label{CTV1} 
\underset{\bold{x},\bold{p}_x,{\bold{p}_y}\in \mathcal{X},\bold{p}_z}{\min}~
 &\sum_{i=1}^N  \sqrt{|\bold{p}_x|_i^2 + |\bold{p}_z|_i^2} + \frac{\lambda}{2} \|\bold{G}\bold{x}-\bold{y}\|_2^2  \\
&\text{subject to}~~
\begin{cases}
\bold{p}_x=\nabla_{\!x} \bold{x}, \nonumber\\
 \bold{p}_y={\bold{x}}, \nonumber \\ 
\bold{p}_z=\nabla_{\!z} \bold{x}. \nonumber
\end{cases}
\end{align}

2) \textit{Augmented Lagrangian.} The second trick is to relax these new linking constraints with an augmented Lagrangian function. This recasts \ref{CTV1} as the following min-max optimization problem
\begin{align}  \label{CTV2}
\underset{\bold{x},\bold{p}_x,\bold{p}_y\in \mathcal{X},\bold{p}_z}{\min}~&\underset{\bold{q}_x,\bold{q}_y,\bold{q}_z}{\max}
 \sum_{i=1}^N  \sqrt{|\bold{p}_x|_i^2 + |\bold{p}_z|_i^2} + \frac{\lambda}{2} \|\bold{G}\bold{x}-\bold{y}\|_2^2 \nonumber \\
&+ \langle \bold{q}_x,\bold{p}_x -\nabla_{\!x} \bold{x}\rangle +   \frac{\xi}{2}\|\bold{p}_x -\nabla_{\!x} \bold{x}\|_2^2 \nonumber \\
&+ \langle \bold{q}_y,\bold{p}_y-\bold{x}\rangle +  \frac{\xi}{2} \|\bold{p}_y-\bold{x}\|_2^2 \nonumber \\
&+ \langle \bold{q}_z,\bold{p}_z -\nabla_{\!z} \bold{x}\rangle +   \frac{\xi}{2}\|\bold{p}_z -\nabla_{\!z} \bold{x}\|_2^2, 
\end{align}
where $\bold{q}_x, \bold{q}_y, \bold{q}_z$ are Lagrange multipliers (dual variables), and $\xi>0$ is the penalty parameter. Note that, here we used the same penalty parameter for all three constraints but one may use different parameter for each of them.
The augmented Lagrangian method \citep[a.k.a. method of multipliers,][]{HESTENES_1969_MAG} maximizes the objective in equation \ref{CTV2} with respect to the dual variables iteratively by using  a simple steepest ascent algorithm (with step length $\xi$)
\begin{align} \label{duals}
\bold{q}_x^{k+1}&=\bold{q}_x^{k} + \xi (\bold{p}_x^{k+1} -\nabla_{\!x} \bold{x}^{k+1})\\
\bold{q}_y^{k+1}&=\bold{q}_y^{k} + \xi (\bold{p}_y^{k+1}-\bold{x}^{k+1}) \\
\bold{q}_z^{k+1}&=\bold{q}_z^{k} + \xi (\bold{p}_z^{k+1} -\nabla_{\!z} \bold{x}^{k+1})
\end{align}
where $\bold{x}^{k+1}$ and $\bold{p}_x^{k+1}, \bold{p}_y^{k+1}, \bold{p}_z^{k+1}$ are obtained by solving 
\begin{align}  \label{CTV2_mm}
\underset{\bold{x},\bold{p}_x,\bold{p}_y\in \mathcal{X},\bold{p}_z}{\arg\min}~&
 \sum_{i=1}^N  \sqrt{|\bold{p}_x|_i^2 + |\bold{p}_z|_i^2} + \frac{\lambda}{2} \|\bold{G}\bold{x}-\bold{y}\|_2^2 \nonumber \\
&+ \langle \bold{q}_x^{k},\bold{p}_x -\nabla_{\!x} \bold{x}\rangle +   \frac{\xi}{2}\|\bold{p}_x -\nabla_{\!x} \bold{x}\|_2^2 \nonumber \\
&+ \langle \bold{q}_y^{k},\bold{p}_y-\bold{x}\rangle +  \frac{\xi}{2} \|\bold{p}_y-\bold{x}\|_2^2 \nonumber\\
&+ \langle \bold{q}_z^{k},\bold{p}_z -\nabla_{\!z} \bold{x}\rangle +   \frac{\xi}{2}\|\bold{p}_z -\nabla_{\!z} \bold{x}\|_2^2,
\end{align}
beginning with $\bold{p}_x^0=\bold{p}_y^0= \bold{p}_z^0=\bold{0}$.

Equations \ref{duals}-\ref{CTV2_mm} can be simplified by using a change of variables ($\bold{q}_x\leftarrow\frac{1}{\xi}\bold{q}_x, \bold{q}_y\leftarrow\frac{1}{\xi}\bold{q}_y, \bold{q}_z\leftarrow\frac{1}{\xi}\bold{q}_z$) using the fact that for two real vectors $\bold{a}$ and $\bold{b}$ the following holds: 
\begin{equation}
\langle \bold{a},\bold{b} \rangle + \frac{\xi}{2} \|\bold{b}\|_2^2= \frac{\xi}{2} \|\bold{b}+\frac{1}{\xi}\bold{a}\|_2^2-\frac{\xi}{2} \|\frac{1}{\xi}\bold{a}\|_2^2.
\end{equation}
Accordingly, 
\begin{align} \label{duals2}
\bold{q}_x^{k+1}&=\bold{q}_x^{k} + \bold{p}_x^{k+1} -\nabla_{\!x} \bold{x}^{k+1}\\
\bold{q}_y^{k+1}&=\bold{q}_y^{k} + \bold{p}_y^{k+1}-\bold{x}^{k+1} \\
\bold{q}_z^{k+1}&=\bold{q}_z^{k} + \bold{p}_z^{k+1} -\nabla_{\!z} \bold{x}^{k+1}
\end{align}
and 
\begin{align}  \label{CTV2_mm2}
\underset{\bold{x},\bold{p}_x,\bold{p}_y\in \mathcal{X},\bold{p}_z}{\arg\min}~&
 \sum_{i=1}^N  \sqrt{|\bold{p}_x|_i^2 + |\bold{p}_z|_i^2} + \frac{\lambda}{2} \|\bold{G}\bold{x}-\bold{y}\|_2^2 \nonumber \\
&+ \frac{\xi}{2}\|\bold{p}_x -\nabla_{\!x} \bold{x} + \bold{q}_x^{k}\|_2^2-\frac{\xi}{2} \|\bold{q}_x^{k}\|_2^2 \nonumber \\
&+ \frac{\xi}{2} \|\bold{p}_y-\bold{x} + \bold{q}_y^{k}\|_2^2-\frac{\xi}{2} \|\bold{q}_y^{k}\|_2^2 \nonumber\\
&+ \frac{\xi}{2}\|\bold{p}_z -\nabla_{\!z} \bold{x} + \bold{q}_z^{k}\|_2^2-\frac{\xi}{2} \|\bold{q}_z^{k}\|_2^2.
\end{align}

3) \textit{Alternating minimization.}
The basic augmented Lagrangian method minimizes the objective function in equation \ref{CTV2_mm2} (augmented Lagrangian function) jointly over $\bold{x}, \bold{p}_x, \bold{p}_y$, and $\bold{p}_z$, the third trick is to perform this minimization by alternating minimizing with respect to each variable separately \citep{Goldstein_2009_SBM,Boyd_2011_DOS} to arrive at the so-called ADMM.
\begin{align}  \label{x-sub}
\bold{x}^{k+1} &=\underset{\bold{x}}{\arg\min}~~
 \frac{\lambda}{2} \|\bold{G}\bold{x}-\bold{y}\|_2^2+ \frac{\xi}{2}\|\bold{p}^k_x -\nabla_{\!x} \bold{x} + \bold{q}_x^{k}\|_2^2 \nonumber \\
& + \frac{\xi}{2} \|\bold{p}^k_y-\bold{x} + \bold{q}_y^{k}\|_2^2 + \frac{\xi}{2}\|\bold{p}^k_z -\nabla_{\!z} \bold{x} + \bold{q}_z^{k}\|_2^2.
\end{align}
\begin{align}  \label{pxpz_sub}
(\bold{p}^{k+1}_x,\bold{p}^{k+1}_z)&=\underset{\bold{p}_x,\bold{p}_z}{\arg\min}
 \sum_{i=1}^N  \sqrt{|\bold{p}_x|_i^2 + |\bold{p}_z|_i^2}  \nonumber \\
&+ \frac{\xi}{2}\|\bold{p}_x -\nabla_{\!x} \bold{x}^{k+1} + \bold{q}_x^{k}\|_2^2 \nonumber \\
&+ \frac{\xi}{2}\|\bold{p}_z -\nabla_{\!z} \bold{x}^{k+1} + \bold{q}_z^{k}\|_2^2.
\end{align}
\begin{align}  \label{py_sub}
\bold{p}_y^{k+1}=\underset{\bold{p}_y\in \mathcal{X}}{\arg\min} ~~\frac{\xi}{2} \|\bold{p}_y-\bold{x}^{k+1} + \bold{q}_y^{k}\|_2^2.
\end{align}
Subproblem \ref{x-sub} is an easy-to-solve least-squares problem, which has a closed-form solution obtained by setting the derivative of the objective function with respect to $\bold{x}$ equal to zero.
\begin{eqnarray}
 \label{CTV4}
&&\bold{x}^{k+1} = 
\big[
\lambda\bold{G}^T \bold{G} + \xi \nabla_{\!x}^T\nabla_{\!x}+\xi \bold{I} +\xi \nabla_{\!z}^T\nabla_{\!z}
\big]^{-1} \\
&&\big[\lambda\bold{G}^T \bold{y} + \xi \nabla_{\!x}^T [\bold{p}_x^k+\bold{q}_x^k]+\xi [\bold{p}_y^k+\bold{q}_y^k]+\xi \nabla_{\!z}^T [\bold{p}_z^k+\bold{q}_z^k]\big]. \nonumber
\end{eqnarray}
Subproblem  \ref{pxpz_sub} also has a closed form solution, given by the generalized soft thresholding function \citep{Goldstein_2009_SBM}, 
 \begin{equation}
\label{px}
\bold{p}_x^{k+1}= \max(1 - \frac{1/\xi}{\bold{r}},0)\circ (\nabla_{\!x} \bold{x}^{k+1}-\bold{q}_x^k),
\end{equation} 
and
\begin{equation}
\label{pz}
\bold{p}_z^{k+1}= \max(1 - \frac{1/\xi}{\bold{r}},0)\circ (\nabla_{\!z} \bold{x}^{k+1}-\bold{q}_z^k),
\end{equation} 
where
\begin{equation} \label{r}
\bold{r}=\sqrt{|\nabla_{\!x} \bold{x}^{k+1}-\bold{q}_x^k|^2+|\nabla_{\!z} \bold{x}^{k+1}-\bold{q}_z^k|^2}.
\end{equation}
Subproblem \ref{py_sub} is a projection operator given by
  \begin{equation} \label{projp}
 \bold{q}^{k+1} = \text{proj}_{\mathcal{X}} (\bold{x}^{k+1} - \bold{q}_y^k).
 \end{equation} 
In the case that $\mathcal{X}$ is a box set of form
\begin{equation}
\mathcal{X}=\{\bold{x}\vert \bold{x}_{min} \leq \bold{x}\leq \bold{x}_{max} \} ,
\end{equation}
then the projection operator admit a closed form solution 
\begin{equation} \label{projp_box}
\text{proj}_{\mathcal{X}} (\bold{x}^{k+1} - \bold{q}_y^k)=\min(\max(\bold{x}^{k+1} - \bold{q}_y^k,\bold{x}_{min}), \bold{x}_{max}),
\nonumber
\end{equation} 
 where $\bold{x}_{min}$ and $\bold{x}_{max}$ are lower and upper bounds of $\bold{x}$, respectively.  \\

%

\bibliographystyle{seg}
\newcommand{\SortNoop}[1]{}

\end{document}